\documentclass[10pt]{article}
\usepackage{amsmath, amssymb}
\usepackage{latexsym}
\usepackage{mathtools}
\usepackage{ulem}
\usepackage{multirow}
\usepackage{graphicx}
\usepackage{color}
\usepackage{url}
\usepackage{framed}
\usepackage{float}
\usepackage{enumitem}
\usepackage{lipsum}

\providecommand{\keywords}[1]{\textbf{Keywords:} #1}

\numberwithin{equation}{section}
\DeclareMathAlphabet{\itbf}{OML}{cmm}{b}{it}

\newcommand{\RR}{\mathbb{R}}

\newcommand{\NN}{\mathbb{N}}
\newcommand{\ZZ}{\mathbb{Z}}

\newcommand{\ds}{\displaystyle}
\newcommand{\no}{\nonumber}
\newcommand{\ri}{\rightarrow}

\newcommand{\bm}{{\itbf m}}

\newcommand{\bx}{{\itbf x}}

\newcommand{\bg}{{\itbf g}}
\newcommand{\bh}{{\itbf h}}
\newcommand{\bev}{{\itbf e}}

\newcommand{\bu}{{\itbf u}}
\newcommand{\bn}{{\itbf n}}
\newcommand{\by}{{\itbf y}}

\newcommand{\bi}{\begin{itemize}}
\newcommand{\ei}{\end{itemize}}

\newcommand{\cu}{{\cal U}}
\newcommand{\cv}{{\cal V}}
\newcommand{\cg}{{\cal G}}
\newcommand{\cf}{{\cal F}}
\newcommand{\ca}{{\cal A}}
\newcommand{\ch}{{\cal H}}
\newcommand{\ce}{{\cal E}}

\newcommand{\bH}{{\itbf H}}

\newcommand{\be}{\begin{eqnarray}}
\newcommand{\ee}{\end{eqnarray}}

\newcommand{\ben}{\begin{eqnarray*}}
\newcommand{\een}{\end{eqnarray*}}
\def\ds{\displaystyle}

\newcommand\ov{\overline}

\newtheorem{lem}{Lemma}[section]
\newtheorem{prop}{Proposition}[section]

\newtheorem{thm}{Theorem}[section]

\newcommand{\bea}{\begin{eqnarray*}}
\newcommand{\eea}{\end{eqnarray*}}
\newcommand{\bean}{\begin{eqnarray}}
\newcommand{\eean}{\end{eqnarray}}
\newcommand{\p}{\partial}
\newcommand{\f}{\frac}

\newcommand{\di}{\mbox{div }}
\newcommand{\aaa}{\mbox{$[$}}
\newcommand{\bbb}{\mbox{$]$}}

\begin{document}

% in this version the tangential vector field  are more 
% explicitly defined
\title{A stochastic algorithm for fault inverse problems in elastic   half space 
with proof of convergence}

\author{  Darko
Volkov \thanks{Department of Mathematical Sciences,
Worcester Polytechnic Institute, Worcester, MA 01609.
}  }
%joancs@cttc.upc.edu
%Heat and Mass Technological Center (CTTC). Technical University
%of Catalonia (UPC), Colom 11, 08222 Terrassa (Barcelona), Spain

\maketitle

\begin{abstract}
A general stochastic algorithm for solving mixed linear and nonlinear problems
was introduced in \cite{volkov2020stochastic}.
We show in this paper how it can be used to solve the fault inverse problem,
where a planar fault in elastic half-space and a slip on that fault have to be reconstructed
from noisy surface displacement measurements. 
With  the parameter giving the plane containing the fault
denoted by $\bm$ and the regularization parameter for the linear part of the inverse problem
denoted by $C$, 
both modeled as random variables, 
we derive a formula for 
the  posterior marginal of $\bm$.
Modeling $C$
 as a random variable allows 
to sweep through a wide range of possible values which was shown
to be superior to selecting  a fixed value \cite{volkov2020stochastic}.
We prove that this posterior marginal of $\bm$ is convergent as the number of 
 measurement points  and the dimension of the space for discretizing slips 
increase. Simply put, our proof only assumes that 
the regularized discrete error functional for processing measurements
relates to an order 1 quadrature rule and that the union of 
the  finite-dimensional
spaces for discretizing slips  
is dense.  
Our proof relies on trace class operator theory to show that
an adequate sequence of determinants is uniformly bounded.
We also explain how our proof can be extended to 
a whole class of inverse problems, as long as some basic requirements are met.
Finally, we  show numerical simulations that illustrate the numerical convergence of our algorithm.
\end{abstract}

\bigskip
\section{Introduction}
In \cite{volkov2020stochastic}, we introduced 
a numerical method for mixed linear and nonlinear inverse problems.
This method applies to cases  where  the data for the inverse problem is corrupted by noise 
and where for each value of the nonlinear parameter, the underlying linear problem is
ill-posed. Accordingly,  regularizing this linear part is required. The norm  used 
for the regularization process has to be multiplied by a scaling parameter, 
also called regularization parameter,
denoted by 
$C$ throughout this paper. In \cite{volkov2020stochastic}, a Bayesian  approach
was adopted, and $C$ was  modeled as a random variable. 
%The prior distribution of $C$ was  assumed to be independent of the prior distribution of
%the nonlinear parameter, $\bm$. 
It was shown in  \cite{volkov2020stochastic} that this approach is superior to
selecting $C$ using some standard method for linear inverse problems, such as 
 the discrepancy principle, or the 
generalized cross validation. Loosely speaking, this can be explained by
observing that for different values of the nonlinear parameter $\bm$, these classical methods
will favor different values of $C$, and as a result different values of the nonlinear parameter $\bm$
cannot be fairly compared. 
Attempting to select a unique value of $C$ for all values of $\bm$ leads to somehow better results,
but as demonstrated in the last section of this paper, doing so pales in comparison to the method
advocated in \cite{volkov2020stochastic}.\\
In this paper, we first derive 
in section \ref{Model derivation}
 a specific version of the
Bayesian posterior distribution of $\bm$ following the 
 method introduced in 
  \cite{volkov2020stochastic}. This version applies to an inverse problem in half space for the
	linear elasticity equations. 
	In the direct problem, a slip field on an open surface that we will call a fault, 
	gives rise to a displacement field. In the inverse  problem, this field 
	is measured on the plane on  top of the half space at $M$  points.
	%The reason for adopting this notation is that further in this paper when
	%we turn to analyzing the convergence of the reconstruction algorithm,
	%$M_N$ 
	%will be a sequence indexed by $N$ and will tend to infinity. 
	%The set of $M_N$ points will not be required to be included in the set of $M_{N+1}$
	%points; instead, these sets of points will  be required to relate to special quadrature rules.\\
	The linear part of the inverse problem consists of reconstructing the slip  field on the fault. 
	The nonlinear part consists of finding the geometry and the location of the fault.
	This formulation is commonly used in geophysics to model slow slip events in the vicinity
	of subduction zones, or the total displacement resulting from a dynamic earthquake,
		see \cite{volkov2017reconstruction, volkov2017determining} and references therein.
	In this paper, the geometry of the fault is assumed to be planar, thus we choose
	the nonlinear parameter $\bm$ to be in $\RR^{3}$.
	The stochastic model considered in this paper is different from a model studied
	in an earlier paper by the same author \cite{volkov2019stochastic}. As explained in
	section \ref{Model derivation}, the difference is that we here  assume the 
	regularization parameter $C$ for the linear part of the inverse problem to be a random variable
	and the covariance $\sigma$ of the measurements is unknown. 
	In section \ref{Model derivation}, the likelihood of $\sigma$ is optimized
	based on the data and as a consequence the resulting posterior distribution function
	of $\bm$ is entirely different from the simpler one  used in \cite{volkov2019stochastic}. 
	Computing this new posterior is more intricate since it involves the random variable
	$(\bm, C)$ instead of just $\bm$ in \cite{volkov2019stochastic}. 
The benefit of using this new posterior is that it leads to much better
	numerical results if  the covariance $\sigma$ of the measurements is unknown as
	 shown in sections
 \ref{Numerical solution} and
 \ref{Failure at fixed C}. \\	
	In section \ref{Convergence result},
	we provide a  mathematical proof of the soundness of our Bayesian approach
	for computing the posterior probability density of $\bm$.
	More precisely, we show that as the number of measurement points $M$
	tends to infinity and the dimension of the space %$H$ 
	for discretizing slip fields 
	tends to infinity,  
	%the marginal posterior distribution function
%of $\bm$ converges uniformly to zero away from the true value $\tilde{\bm}$.
%Remarkably, this result holds even at high values for noise covariance.
the probability of $\bm$ to be further than a fixed $\eta>0$
from the true value $\tilde{\bm}$ converges to zero if the noise level is low enough.
Interestingly, although the derivation of the probability law of the posterior of $\bm$
assumes that the noise is Gaussian, once this law is set the proof of convergence does not require the
noise to be Gaussian.
The proof assumes that 
the regularized discrete error functional for sampling measurements
relates to an order 1 quadrature rule for $C^1$ functions
(which is verified by most commonly used quadrature rules)
and that the union of 
the  finite-dimensional
spaces for discretizing slips  
is dense, a natural requirement that can be easily satisfied by combining
adequate finite element spaces. 
The more difficult step in this proof of convergence requires
showing that
a sequence of determinants is uniformly bounded: this can be achieved 
by using trace class operator theory \cite{gohberg1996traces}.\\
In section \ref{extension}, for the reader's convenience, 
we extend the formulation of our recovery method to more general inverse problems.
We write down the expression of the posterior probability distribution that is at the core
of our reconstruction method. A companion 
parallel sampling algorithm is  given in  \cite{volkov2020stochastic}, section 3.2.
At the end of section \ref{extension} of this paper, we
state a convergence theorem for the posterior probability distribution
which is valid in this more general framework. \\
% list the necessary assumptions on a general inverse problem for our proof of convergence 
%for  the posterior of the nonlinear parameter to hold. \\
In section \ref{Numerical simulations}, we present 
numerical simulations that use our posterior probability distribution  for reconstructing the
geometry parameter $\bm$. The posterior marginals of $\bm$ and of the regularization parameter 
$C$ are computed using the method of choice sampling due to the size of the search space.
 We used a modified version of the Metropolis algorithm
which is well suited to parallel computing, \cite{calderhead2014general}.
This modified version was  described in details in previous work, \cite{volkov2020stochastic}.
The computed posterior marginals shown in  section \ref{Numerical simulations}
illustrate the theoretical convergence result proved in section \ref{Convergence result}:
in the low noise case, as well as in the high noise case, 
the posterior marginals of $\bm$ tighten around the  true value
 $\tilde{\bm}$ as   the number of measurement points %$M_N$
and the dimension of the space %$H_p$  
for discretizing slips increase,
as expected from  Theorem \ref{main cv theorem}.
Of course in the low noise case, the tightening is more narrow,   also  as expected from Theorem \ref{main cv theorem}.
In a final section, we show how using a simplistic method where different values 
of the regularization parameter $C$ are first fixed and then  marginal posteriors of $\bm$ 
are computed, leads to results that are impossible to interpret since there is no objective
way to select an optimal $C$ for mixed linear and nonlinear problems if the variance
of the noise is unknown.

\section{Governing equations, inverse problem, and error functionals} \label{num sim}

\subsection{Formulation of the direct and the inverse fault problems}
%We now formulate the
%direct and the inverse fault  problems/
 %half space fault inverse problem in seismology that we will solve
%using the algorithm introduced in section \ref{par} .
Using  standard rectangular coordinates, let $\bx= (x_1, x_2, x_3)$ denote elements of 
 $\RR^{3}$.
%A point $(x,y,z)$ will denote a point 
We define $\RR^{3-}$ to be the open half space $x_3<0$.
%Let  $\Gamma$ be a bounded open surface
%with smooth boundary included in the plane $\{ x_3=0 \}$.
%Let $\p_i$ denote the derivative in the $i$-th coordinate.
%In this paper 
%we  consider the case of linear, homogeneous, isotropic elasticity; 
The direct problem relies on  the equations of linear elasticity with
Lam\'e constants $\lambda$ and $\mu$ 
such that $\lambda > 0$ and $\lambda + \mu> 0$.
For a vector field $\cv = (\cv_1, \cv_2, \cv_3)$,
%the stress  tensors will be denoted,
%\bea
%\sigma_{ij}(\cv) = \lambda \, \di \cv \, \delta_{ij} + \mu \, (\p_i \cv_j + \p_j \cv_i ), \\
%\eea
%and 
the  stress vector in the 
direction $\bev \in \mathbb R^3 $ will be denoted by
\bea
T_\bev \cv = \sum_{j=1}^3
\left(  \lambda \, \di \cv \, \delta_{ij} + \mu \, (\p_i \cv_j + \p_j \cv_i ) \right)
 e_j.
\eea
Let $\Gamma$ be a Lipschitz open surface strictly included in $\RR^{3-}$
with  normal
vector $\bn$ defined almost everywhere.
We define the jump  $\aaa \cv \bbb$ of the  vector field 
$\cv$ across $\Gamma$ 
to be
$$
\aaa \cv \bbb (\bx) = \lim_{h \ri 0^+} \cv  (\bx + h \bn) - \cv  (\bx - h \bn),
$$
for $\bx$ in $\Gamma$, if this limit exists.
Let $\cu$ be the displacement field solving
\bean
\mu \Delta \cu+ (\lambda+\mu) \nabla \di \cu= 0  \mbox{ in } \RR^{3-} 
\setminus \Gamma \label{uj1}, \\
 T_{\bev_3} \cu =0 \mbox{ on the surface } x_3=0 \label{uj2}, \\
 T_{\bn} \cu  \mbox{ is continuous across } \Gamma \label{uj3}, \\
 \aaa  \cu\bbb =\cg \mbox{ is a given jump across } \Gamma , \label{uj3andhalf}\\
\cu (\bx) = O(\f{1}{|\bx|^2}),  \nabla \cu (\bx) = O(\f{1}{|\bx|^3}), \mbox{ uniformly as } 
|\bx| \rightarrow \infty,
\label{uj4}
\eean
where %$\bx=(x_1, x_2, x_3)$ is any point in $\RR^{3-}$ and
$\bev_3$ is the vector $(0,0,1)$.
Let 
  $\widetilde{H}^{\f12}(\Gamma)^2$ 
  be the space  of restrictions to $\Gamma$ of tangential fields  in
   $H^{\f12}(\p D)^2$ supported in $\ov{\Gamma}$, where
  $D$ is  a bounded domain in $\RR^{3-}$ 
such that $\Gamma \subset \p D$. 
In \cite{volkov2017reconstruction}, we defined the functional space 
${\bf S}$ of vector fields $\cv $ defined in $\RR^{3-}\setminus{\ov{\Gamma}}$
such that $\nabla \cv $ and $\ds \f{\cv}{(1+ r^2)^{\f12}}$ 
are in $L^2(\RR^{3-}\setminus{\ov{\Gamma}})^3$ and we proved the following
	result.
\begin{thm}\label{direct exist unique}
 Let $\cg$ be in $\widetilde{H}^{\f12}(\Gamma)^2$.
The problem
(\ref{uj1}-\ref{uj3andhalf})
 has a unique solution in ${\bf S}$.
In addition, the solution $\cu$ satisfies the decay conditions
(\ref{uj4}).
\end{thm}
%We recall that if $\cg$ is continuous, the support of $\cg$, 
%$\mbox{supp } \cg$, is equal to the closure
%of the set of points in $\Gamma$ where $\cg$ is non zero; in general $\mbox{supp }  \bg$
%is defined in the sense of distributions. 
%Can  both $\cg$ and $\Gamma$  be determined from  given only on the plane $x_3=0$?
The following theorem shown in \cite{volkov2017reconstruction}
asserts that
 $\cg$ and $\Gamma$   are uniquely determined 
from the data
 $\cu$
 given  on a relatively open
set  of the plane $x_3=0$ if we know that $\Gamma$ is planar.
\begin{thm}{\label{uniq1}}
Let $\Gamma_1$ and $\Gamma_2$ be two 
planar  open 
surfaces.
%, with smooth boundary, such that each of them  is included in a
%rectangle strictly contained in $\RR^{3-}$.
For $i$ in $\{ 1,  2\}$, assume that $\cu^i$ solves  (\ref{uj1}-\ref{uj4}) for 
$\Gamma_i$ in place of $\Gamma$ and $\cg^i$, a tangential field in 
$\tilde{H}^{\f12}(\Gamma_i)^2$, in place of $\cg$.
Assume that $\cg^i$ has full support in $\Gamma_i$, that is, 
$\mbox{supp } \cg_i = \ov{\Gamma_i}$.
Let $V$ be a non empty open subset in $\{x_3 =0\}$.
If $\cu^1$ and $\cu^2$ are equal in $V$, then
$\Gamma_1 =  \Gamma_2$
and $\cg^1=\cg^2$.
\end{thm}
Theorems  \ref{direct exist unique} and \ref{uniq1} 
were  proved in \cite{volkov2017reconstruction} 
for media with constant Lam\'e coefficients
and in \cite{aspri2020analysis} for more general Lam\'e
systems.
Similar results were obtained in \cite{aspri2020dislocations} for layered media, albeit
in a bounded domain.\\
%. Later,
%in \cite{aspri2020analysis}, 
%the direct problem (\ref{uj1}-\ref{uj3andhalf})
% was analyzed under weaker conditions 
%for the regularity of the displacement field $\cu$ and the slip
%field $\cg$. In \cite{aspri2020arxiv},
%the direct problem (\ref{uj1}-\ref{uj3andhalf})
%was proved to be uniquely solvable in case of
%piecewise Lipschitz coefficients and general elasticity tensors.
%Both \cite{aspri2020analysis} and \cite{aspri2020arxiv}
%include a proof of  uniqueness  for the  fault inverse problem under appropriate assumptions.
%In our case, 
The  solution $\cu$ to problem (\ref{uj1}-\ref{uj3andhalf})  can be expressed  
 as the convolution on $\Gamma$
\bean
  \cu(\bx) = \int_\Gamma \bH(\bx, \by, \bn) \cg(\by) \, d s (\by) \label{int formula},
\eean
where $\bH$ is  the  Green tensor  associated to the system (\ref{uj1}-\ref{uj4}), and
  $\bn$ is the vector normal to $\Gamma$. 
The practical determination of this half space Green tensor $\bH$ was  studied
in  \cite{Okada} and later, more rigorously, in \cite{volkov2009double}.
%In particular, $\bH$ satisfies the decay conditions
%\bea
% \bH(\bx, \by, \bn) = O(|\bx|^{-2}), \, \nabla_x  \bH(\bx, \by, \bn) = O(|\bx|^{-3}), \q |\bx| \ri \infty,
%\eea
%uniformly in $\by$ and in $\bn$, as long as $\by$ remains in a bounded subset of
%$\RR^{3-}$.
Due to formula (\ref{int formula}) %and the  decay of the Green's tensor
%$H$ at infinity, 
we can define a continuous mapping  from tangential fields
$\cg $ in $H^{1}_0(\Gamma)^2$ to surface displacement fields $\cu(x_1, x_2, 0)$
in $L^2(V)^3 $ where $\cu$ and $\cg$ are related
by (\ref{uj1}-\ref{uj4}). This mapping is  compact since $\bH(\bx, \by, \bn)$ is smooth for
$(\bx, \by)$ in $\ov{V}\times \ov{\Gamma}$.
Theorem \ref{uniq1} asserts that 
this compact mapping is injective, so its inverse  can be defined.
As the inverse of a compact mapping is unbounded,
 % \cap C^\infty (\RR^2)$. 
%It is well known, however, that such an operator ${\cal M}$ is compact, therefore
%its inverse is unbounded. 
%It is thus clear that a stable numerical method for reconstructing
%$\cg$ from $\cu(x_1, x_2, 0)$ will have to use some regularization process.
finding $\cg$ from $\cu(x_1, x_2, 0)$ has to involve regularization.
%This clearly points to the need to use regularization to find $\cg$ from $\cu(x_1, x_2, 0)$.

\subsection{A regularized functional for the   reconstruction of  planar faults}
\label{Continuous regularized reconstruction} 
Let $R$ be a bounded, non-empty, open set  of the plane $x_3=0$.
Let $B$ be a compact set of $(a, b, d)$ in $\RR^3$ such that 
\bea 
 \{  (x_1, x_2, a x_1 + b x_2 +d): (x_1, x_2) \in R\}
\eea 
is included in the half-space $x_3 <0$.
We introduce the notations
\bea
\bm=(a,b,d), \\
\Gamma_\bm =\{  (x_1, x_2, a x_1 + b x_2 +d): (x_1, x_2) \in R\}.
\eea
%We assume that $B$ is a closed and bounded subset of $\RR^3$.
 It follows  
that 
%\bean \label{pos dis}\begin{array}{l}\mbox{\sl
 the distance between $\Gamma_{\bm}$ and the plane $x_3=0$ is bounded below
%}\\\mbox{\sl 
by the same  positive constant for all $ \bm$ in $B$.
%}\end{array}\eean
We assume that slips $\cg$ are supported in such sets $\Gamma_{\bm}$.
% (meaning that their supports are included in 
%$\Gamma_m$, but they could be different from $\Gamma_\bm$).
We can then map all these fields into the rectangle $R$.
We thus obtain displacement vectors for $\bx$ in $V$ by the integral formula
\bean \label{new int}
\cu(\bx, \cg, \bm) = \int_R \bH_\bm(\bx, y_1, y_2)
 \cg (y_1, y_2) s d y_1 d y_2, 
\eean
for any $\cg$ in $H^1_0 (R)^2$ and  $\bm$ in $B$, where $s $ is the surface element on $\Gamma_\bm$
and  $\bH_\bm(\bx, y_1, y_2)$ is derived from 
the Green's tensor $\bH$ for $\by$ on $\Gamma_\bm$.
%Let $V$ be a non-empty, open subset of the plane
%$x_3 =0 $ and 
We now assume that $V$ is a bounded open subset of the plane $x_3 =0$.  For 
a fixed $\tilde{\cu}$ be in $L^2(V)$,
and  a fixed $\bm$ in $B$ we define the regularized error functional
\bean
 \cf_{\bm,C} (\cg)= \int_V  |\cu(\bx, \cg, \bm) - \tilde{\cu}(\bx ))|^2
 d\bx
+ C \int_{R} |\nabla \cg|^2 ,
\label{F cont abd}
\eean
where
% ${\cal C} (\bx) $ is a diagonal positive definite 3 by 3 matrix 
%for $\bx$ in $\ov{V}$, which is continuous in $\bx$, 
% and 
$C>0$ is the regularization parameter and $\cg$ is in $H^1_0(R)^2$. 
%In formula (\ref{F cont abd}) we intentionally
%used ${\cal C}^{-1}$ rather than ${\cal C}$ because we will later view it as a covariance term.
 Define the operator
\bean 
\ca_{\bm} &:& H^1_0 (R)^2 \ri L^2(V)^3 \no \\
&&\cg \ri \int_R  \bH_\bm(\bx, y_1, y_2)
 \cg (y_1, y_2) s d y_1 dy_2 .  \label{Aabd}
\eean
It is clear that $\ca_{\bm} $ is linear, continuous, and compact.
The functional $\cf_{\bm,C}$ can also be written as,
\bean
 \cf_{\bm,C} (\cg)= \|\ca_{\bm} \cg  - \tilde{\cu} \|^2_{L^2(V)^3} + C \| \cg \|^2_{H^1_0(R)^2},
\eean
where
% in $L^2(V)$ we use the norm
%\bean \label{normV}
%\| \bu \|_{L^2(V)} = 
%(\int_V  \bu(\bx )'
%{\cal C}^{-\f12}(\bx)  \bu(\bx) 
 % d\bx)^{\f12},
%\eean
 in $H^1_0(R)^2$ we use the norm
\bean \label{normR}
\| \cg \|_{H^1_0(R)^2} =
(\int_{R} |\nabla \cg|^2 )^{\f12}.
\eean
In the remainder of this paper, for the sake of simplifying notations,
both $\|   \, \|_{L^2(V)^3}$ and $\| \,  \|_{H^1_0(R)^2} $ will be abbreviated by $\| \, \|$;
context will eliminate any risk of confusion. 
\begin{prop}\label{Fmin}
For any fixed $\bm$ in $B$ and $C>0$, $\cf_{\bm,C} $ achieves a unique minimum $\ch_{\bm,C}$
in $H^1_0 (R)^2$.
\end{prop}
\textbf{Proof}:
This result holds thanks to classic Tikhonov regularization theory
(for example, see \cite{kress1989linear}, Theorem 16.4). $\Box$

%On the subject of the widely used Tikhonov regularization technique, we recommend
%the textbook \cite{kress1989linear}, Chapter 16, for  a presentation 
%which is particularly  relevant to our study. 
For $\ch_{\bm,C}$  as in the statement of Proposition \ref{Fmin} we set,
	\bean \label{habdc def}
	 f_C(\bm) =  \cf_{\bm,C} (\ch_{\bm,C}).
	\eean
	\begin{prop} \label{fmin}
	$f_C$ is a Lipschitz continuous function on $B$. It achieves its minimum value on
	$B$.
\end{prop}
\textbf{Proof}:
This was proved in \cite{volkov2019stochastic}, Proposition 3.2. $\Box$\\
It is helpful to state a direct consequence of  uniqueness Theorem \ref{uniq1} in terms of 
the operators $\ca_\bm$.
\begin{prop} \label{injectivity of Am}
Let $\bm, \bm'$ be in $B$ and $\cg, \cg'$ be in $H^1_0(R)^2$.
Assume that $\cg$ or $\cg'$ is non-zero. If $\ca_{\bm}  \cg = \ca_{\bm'}  \cg' $
then $\bm = \bm'$ and $\cg =  \cg'$.
\end{prop}

	In the remainder of this paper,  we only consider one-directional fields $\cg$.
	Accordingly, $\cg$ can be considered to be a scalar function in the space
	$H^1_0(R)$, and  $\ca_\bm$ becomes a linear operator
	from $H^1_0(R)$ to $L^2(V)^3$.

\subsection{A functional for the reconstruction of planar faults from surface measurements
at $M$ points} \label{finite recon}

For $j=1, .. ,M$, let $ P_j$ be points in $V$, so on the surface $x_3 =0 $, and
$\tilde{\bu}(P_j)$ be measured displacements at these points.
%This notation for the measurement points will prove useful for the convergence analysis 
%shown later in this paper; in particular, the set 
%$\{ P_j^N: 1 \leq j  \leq M_N \}$ is not necessarily included in
%$\{ P_j^{N+1}: 1 \leq j  \leq M_{N+1} \}$.
%The number of points $M_N$ is assumed to be increasing in $N \in \NN$ and
%$\lim_{N \ri \infty} M_N = \infty$.
%Let $H_p$ be a  sequence of finite-dimensional subspaces of
%$H_0^1(R)$ such that  $H_p \subset H_{p+1}$ and $\bigcup_{p=1}^\infty H_p$ is dense
%in $H_0^1(R)$.
Let $H$ be a  finite-dimensional subspaces of
$H_0^1(R)$. 
For $\cg$ in  $H$ and $\bm$ in $B$, define the functional 
\bean
F_{\bm,C} (\cg) = \sum_{j=1}^{M} C'(j) | 
(\ca_{\bm} \cg- \tilde{\bu} ) (P_j)|^2 
+ C \int_R |\nabla \cg|^2 ,
\label{Ffunc}
\eean
where 
%$C'=((\sum_{j=1}^{N}  | \tilde{\bu} (P_j^N)  |^2))^{-1}$, $$
$\ca_{\bm} $ was defined in (\ref{Aabd}), and 
$C>0$ is the regularization parameter. 
We assume that the constants $C'(j) $ are positive.
% and 
%form a convergent quadrature rule in $V$, more precisely,
%simply put, they relate 
%$F_{\bm,C} $ to $\cf_{\bm,C} $ as the number of points $M_N$ tends to infinity. 
%More precisely we assume that
% $\cal C$ is smooth and that 
%for all positive integer $k$,  and 
%for all $\varphi$ in $C(\overline{V})$, 
%there is a constant 
%$C(k)$ such that
%\bean \label{quad rule}
%\lim_{N \ri \infty}  \sum_{j=1}^{M_N} C'(j,N) \varphi(P_j^N)
%= \int_V \varphi.
%|\int_V \varphi - \sum_{j=1}^{M_N} C'(j,N) \varphi(P_j^N) |
%\leq C(k ) N^{-\beta} \sup_{V} \sum_{|l| \leq k}|D^l \varphi| ,
%\eean
%where $D^l \varphi $ is a  partial derivative of $\varphi$ with total order $l$
%and $\beta$ is a positive integer depending on $k$.
%We also assume that $C'(j,N) >0$ for all positive integer $N$ and all 
%$j=1, .., N$.
 
	\begin{prop} \label{prop min}
The functional $F_{\bm,C}$ achieves a unique minimum on $H$.
\end{prop}
%	We don't know whether this minimum value is achieved at only one point.
%	Actually, for the value $N=1$ (only one measurement point), this is in general not true.
%	However, as $N$ grows large, we can state an interesting convergence result.
	\textbf{Proof}:
	This results again from Tikhonov regularization theory. $\Box$\\
%(see \cite{kress1989linear}, Theorem 16.4).\\\\
As $F_{\bm,C}$ achieves its minimum
	at some $\bh_{\bm,C}$ in $H$,% (the superscript "disc" stands for discrete).
	we set 
	\bean
	 f^{disc}_C(\bm) =  F_{\bm,C} (\bh_{\bm,C}) \label{define bhdiscmc}.
	\eean
	The superscript "disc" stands for discrete.
	\begin{prop}
	$f^{disc}_C$ is a Lipschitz continuous function on $B$ and achieves its minimum value on
	$B$.
\end{prop}

\textbf{Proof}:
The proof is similar to that of Proposition \ref{fmin}.  $\Box$

%\section{Stochastic model} \label{Stochastic modeling}%247

\section{Bayesian model derivation}\label{Model derivation}
We assume in our stochastic model that the geometry parameter $\bm=(a,b,d)$ is in $B$
%,  the 
% slip field $\bg$ in $H_p$, 
and the measurements $\tilde{\bu} (P_j)$ are related by the equation
\bean \label{noise model}
 %m = \argmin_{ m \in B} \argmin_{\bg \in H^1_0(R)} F^{disc}_{m,C} (\bg)
(\tilde{\bu} (P_1), ..., \tilde{\bu} (P_{M}) ) =  \delta+ 
(\ca_\bm \bg (P_{1}), ..., \ca_\bm \bg (P_{M}) )+ {\cal E},
\eean
where $\ca_\bm$ is as in (\ref{Aabd}), 
$\bg$ is in $H $, 
%%\bm$ and $\bg$ are  now  random variables, 
and ${\cal E}$ in
$\RR^{3 M}$ is additive noise, which 
is modeled to be a random variable. 
The error term $\delta  $ is given by
\bean
\label{delta def}
(\ca_\bm (I- \Pi)(\cg) (P_1), ..., \ca_\bm (I- \Pi)(\cg) (P_M) ),
\eean
where $\Pi$ is the orthogonal projection from 
$H^1_0(R) $ to $H$ and $\cg$ is in $H^1_0(R) $.  $\delta  $ represents the part of  the slip
that cannot be reached by the  finite dimensional  model
and $\bg= \Pi \cg$.
 We assume in this section  that ${\cal E}$ follows a normal  probability distribution
%$\rho_{noise}$ 
with mean zero and diagonal covariance matrix  
$\sigma^2 I$  after rescaling by  $C'$, so altogether the
diagonal of this covariance matrix is the vector in $\RR^{3M}$ is
\bea
\sigma^2 (C'(1)^{-1},C'(1)^{-1},C'(1)^{-1}, ..., 
C'(M)^{-1},C'(M)^{-1},C'(M)^{-1}).
\eea
 %that for a measurement
 %\bean  \label{rhonoise}
%\rho_{noise}(\bv_1, ..., \bv_N) \propto \exp ( - \f12
 %\sum_{j=1}^{N} C'(j,N) | {\cal C}^{-\f12}
%\bv_j|^2 )
%\eean 
The $3 M$ dimensional vector $(\tilde{\bu} (P_1), ..., \tilde{\bu} (P_{M}) ) $
will be denoted by $\tilde{\bu}$. 
Given our assumption on the noise ${\cal E}$, for a fixed $\bm$ and $\bg$ the random variable
$(\ca_\bm \bg (P_{1}), ..., \ca_\bm \bg (P_{M}) )+ {\cal E}$
is normal with mean $(\ca_\bm \bg (P_{1}), ..., \ca_\bm \bg (P_{M}) ) $
and covariance $\sigma^2 I$.
Accordingly, by \eqref{noise model},
the  probability density of 
$\tilde{\bu} - \delta$ knowing the geometry parameter $\bm$, the slip field $\bg$ 
and the variance $\sigma$
is 
\bean \label{umeas}
\rho(\tilde{\bu} - \delta | \bm, \bg, \sigma ) \propto  
\exp ( -\f{1}{2 \sigma^2} \sum_{j=1}^{M} C'(j) 
| 
\ca_\bm \bg- \tilde{\bu}|^2 (P_j) ).
\eean
Next, we assume that  $\bm$ in $B$ and $\bg$ in 
$H$  are independent random variables.
The prior distribution of $\bm$, $\rho_{pr}$, is assumed to be uninformative,
that is, $\rho_{pr}(\bm) \propto 1_B (\bm)$.  
For the prior distribution of  $\bg$ knowing $C$
we follow the Maximum Likelihood (ML) model introduced 
in \cite{galatsanos1992methods}: this is due to the fact that we consider 
$\sigma^2$ to be unknown. In \cite{galatsanos1992methods},
Galatsanos and Katsaggelos only studied a linear problem 
and later 
in  \cite{volkov2020stochastic},
their method was extended to the case where the linear operator depends on  an
unknown nonlinear parameter. In that case it is advantageous to
model the regularization parameter $C$ as a random variable
\cite{volkov2020stochastic}, and the prior of $\bg$ is set to be,
% is 
%assumed to be Gaussian with mean zero and given by
\bean \label{gprior}
%\rho_{H_p}(\bg | \sigma, C) \propto \exp( -\f{C}{2 \sigma^2} \|\bg \|^2),
\rho(\bg | \sigma, C) \propto \exp( -\f{C}{2 \sigma^2} \|\bg \|^2),
\eean
where $\|\bg \|^2$ is given by the square of  the natural norm in $H_p$, $\int_R |\nabla\bg|^2$.
Modeling $C$ as a random variable makes the stochastic formulation derived in this paper
entirely different from an earlier formulation derived in \cite{volkov2019stochastic}.
Indeed, in \cite{volkov2019stochastic}, section 5.4 describes an algorithm for a uniform selection 
of $C$ for all values of $\bm$. That algorithm was certainly an improvement over applying 
classic selection algorithms for $C$ for linear problems (such as Generalized Cross Validation or the discrepancy principle), separately  for each  value of $\bm$,
 which leads to inconsistent regularization demands 
for different values of $\bm$ and an overall poor numerical method. 
However, the algorithm proposed in \cite{volkov2019stochastic} section 5.4, requires
a good knowledge of the expected value of the  magnitude of error measurements.
Instead, in the model adopted here, the prior of $C$ is free to roam over a large interval, 
and the data is used to compute the expected value of $C$. 
Interestingly, these two different approaches lead to entirely different formulas
for posterior probability distributions. 
The approach adopted here leads to much better numerical results as demonstrated in sections
 \ref{Numerical solution} and
 \ref{Failure at fixed C}. The downside of the approach adopted here is that it adds one dimension
to the stochastic variable whose distribution function is to be computed, and accordingly
the computational challenge is one order of magnitude greater.
This downside can be overcome by using a modified version of the Metropolis algorithm
which is well suited to parallel computing 
 \cite{calderhead2014general, volkov2020stochastic}.\\
%\bean
%\rho( m, \bg| \tilde{\bu}_{meas}  ) \propto 
%\ds \rho(\tilde{\bu}_{meas}  | m , \bg) \rho_{{\cal F}_p}(\bg) \rho_{prior} (m)
%.
%\label{Bayes}
%\eean
Set for $\bg$ in $H_p$,
\bean
F_{\bm, C}^\delta (\bg)=  \sum_{j=1}^{M} C'(j) | 
((\ca_{\bm} \bg) - \tilde{\bu} + \delta) (P_j)|^2 
+ C \int_R |\nabla \bg|^2 ,
\label{Ffunc delta}
\eean
and let $\ca_{\bm}^{M }$ be the finite dimensional linear operator 
from $H$ to $\RR^{3M}$ defined by
\bean
\ca_{\bm}^{M} \bg = (C'(1)^{\f12}  \ca_\bm (\bg)(P_1), ..., C'(M)^{\f12}  \ca_\bm (\bg)(P_{M}))
\label{AmpN def}
\eean

\begin{prop}\label{pdf law}
Assume that $(\tilde{\bu} (P_1), ..., \tilde{\bu} (P_{M}) ) -  \delta$
is not zero.
As a function of $\sigma$, $\rho(\tilde{\bu} | \bm, \bg, \sigma )$, 
the probability density of the measurement 
$\tilde{\bu}$ knowing $\bm, \bg$, and $\sigma$,
achieves a unique maximum at
\bean
\sigma_{max}^2= \f{1}{3M} F_{\bm, C}^\delta(\bg_{min})   ,
\label{sigma}
\eean
where $\bg_{min}$ is the minimizer
of $F_{\bm, C}^\delta$ over $H$.
Fixing $\sigma=\sigma_{max}$, the probability density of $(\bm, C)$ knowing $\tilde{\bu}$
 is then
given, up to a multiplicative constant, by the formula
\bean
\rho ( \bm ,C| \tilde{\bu}) \propto  \det (C^{-1} (\ca_{\bm}^{ M })' \ca_{\bm}^{ M} 
 + I)^{-\f12} 
   F_{\bm, C}^\delta (\bg_{min})^{-\f{3M}{2}} \rho_{pr} (\bm, C ) . \label{final}
\eean
\end{prop}
\textbf{Proof:}
Although a slightly different version of this proof can be found in  \cite{volkov2020stochastic},
 we still include it in this paper for the sake of completion.
Using equation \eqref{noise model} and the probability law of
${\cal E}$, the probability density of $\tilde{\bu}$ knowing $\bg$, $\sigma$,
$\bm$,  and $C$ is, 
since $\tilde{\bu}$ does not depend on $C$,
\bean
\rho(\tilde{\bu} |  \bg, \sigma, \bm, C)
= \rho(\tilde{\bu} |  \bg, \sigma, \bm) \no \\
 = (\f{1}{2 \pi \sigma^2})^{\f{3M}{2}} 
\prod_{j=1}^{M} C'(j)^{\f32}
\exp (- \f{1}{2 \sigma^2} 
\sum_{j=1}^{M} C'(j)
|\tilde{\bu}- \ca_\bm \bg  - \delta|^2(P_j)).
\label{u knowing m g sig}
\eean
Recalling \eqref{gprior},
\bean
\rho(\bg |\sigma, \bm, C) = 
\rho(\bg |\sigma,  C)
=   (\f{C}{2 \pi \sigma^2})^{\f{\dim H}{2}}
\exp(- \f{C}{2 \sigma^2} \| \bg \|^2) \label{g prior} ,
\eean
since this prior is independent of $\bm$. 
The joint distribution of $(\tilde{\bu}, \bg)$ knowing $\sigma, \bm, C$ is related
 to the distribution of $\tilde{\bu}$ knowing $\bg, \sigma, \bm, C$ by
\bean
\rho(\tilde{\bu} , \bg | \sigma, \bm, C) = \rho(\tilde{\bu} | \bg , \sigma, \bm, C) 
(\int \rho(\tilde{\bu} , \bg | \sigma, \bm, C) d \tilde{\bu}). \label{to combine}
\eean
Now,  $\int \rho(\tilde{\bu}, \bg | \sigma, \bm, C) d \tilde{\bu} $ is the 
prior 
 probability distribution
of $\bg$ \cite{kaipio2006statistical},
which we said was given by (\ref{g prior}).
%Following \cite{galatsanos1992methods}, 
%we then maximize in $\sigma$ 
%the probability density of $\bu$ knowing $C, \sigma^2, \bm$. 
Combining (\ref{u knowing m g sig}, \ref{g prior}, \ref{to combine})
we obtain,
\bean
\rho (\tilde{\bu}|\sigma, \bm, C) = \int \rho(\tilde{\bu}, \bg | \sigma, \bm, C) d \bg 
=  (\f{1}{2 \pi \sigma^2})^{\f{\dim H+3 M}{2}} C ^{\f{\dim H}{2}}
\prod_{j=1}^{M} C'(j)^{\f32} \no \\
\int \exp(- \f{C}{2 \sigma^2} \| \bg \|^2  - \f{1}{2 \sigma^2} 
\sum_{j=1}^{M} C'(j)
|\tilde{\bu}- \ca_\bm \bg  - \delta|^2(P_j)      ) d \bg .
\label{inter}
\eean
The latter integral can be computed explicitly  \cite{volkov2019stochastic} to find 
\bean
\int \exp(- \f{C}{2 \sigma^2} \| \bg \|^2  - \f{1}{2 \sigma^2} 
\sum_{j=1}^{M} C'(j)
|\tilde{\bu}- \ca_\bm \bg  - \delta|^2(P_j)       ) d \bg 
\no \\
= \exp ( - \f{1}{2 \sigma^2} F_{\bm, C}^\delta(\bg_{min})   ) 
(\det (\f{1}{2 \pi \sigma^2}(  (\ca_{\bm}^{M })'\ca_{\bm}^{M } + C  I)))^{-\f12},
\label{with det}
\eean
where $\bg_{min}$ minimizes  $F_{\bm, C}^\delta$ over $H$ and
$\ca_{\bm}^{ M }$ was defined in  \eqref{AmpN def}.
The determinant in (\ref{with det}) is of order $\dim H$ so some terms in $\sigma$ in 
(\ref{with det}) and (\ref{inter}) simplify to obtain,
\bean
(\f{1}{2 \pi \sigma^2})^{\f{3M}{2}}  C ^{\f{\dim H}{2}}
\prod_{j=1}^{M} C'(j)^{\f32} 
\exp ( - \f{1}{2 \sigma^2} F_{\bm, C}^\delta(\bg_{min})) 
\no \\
(\det ( (\ca_{\bm}^{M })'\ca_{\bm}^{M } + C  I))^{-\f12},
\label{tomax}
\eean
which we now maximize for $\sigma$ in $(0, \infty)$.
Note that $\bg_{min}$ does not depend on $\sigma$.
As $\sigma$ tends to infinity, the limit of (\ref{tomax}) is clearly zero.
As $\sigma$ tends to zero
the limit 
of (\ref{tomax}) is again zero since $(\tilde{\bu} (P_1), ..., \tilde{\bu} (P_{M}) ) -  \delta$
is not zero which implies that $F_{\bm,C}^\delta (\bg_{min}) >0$.
%, as long as $\tilde{\bu}$ is non-zero, $\| \bg_{min} \| \neq 0$, 
%so the limit 
%of (\ref{tomax}) is again zero.
We then take the derivative of (\ref{tomax}) in $\sigma$ and set it to equal to zero to find the equation
\bea
-3M \sigma^{-3M-1} + \sigma^{-3M}(-2)\sigma^{-3} ( - \f{1}{2 }
 F_{\bm, C}^\delta(\bg_{min}   )    )
=0, 
\eea
thus the value
\bea
\sigma_{max}^2= \f{1}{3M} F_{\bm, C}^\delta(\bg_{min})   
%\label{sigma}
\eea
maximizes the density $\rho (\tilde{\bu}| \sigma, \bm, C) $.
Substituting (\ref{sigma}) in 
(\ref{tomax}) we find
for this particular value of $\sigma^2$,
\bean \label{full law}
\rho ( \tilde{\bu}|\bm,C) = 
(\f{3}{2 \pi e})^{\frac{3M}{2}} \prod_{j=1}^{M} (M C'(j))^{\f32} %\no \\
 (\det ( C^{-1}(\ca_{\bm}^{M })'\ca_{\bm}^{M } +   I))^{-\f12}
 F_{\bm, C}^\delta(\bg_{min})^{-\f{3M}{2}} .
\eean
Since our goal is to reconstruct $\bm$ and $C$ knowing $\tilde{\bu}$ we apply Bayes' law 
\bean \label{bayes2}
\rho (\bm , C| \tilde{\bu}) \propto \rho (\tilde{\bu}|\bm, C) \rho_{pr} (\bm, C),
\eean
to obtain (\ref{final}).    $\Box$\\\\

\section{Convergence result}\label{Convergence result}
We  prove in this section a convergence result as the dimension of $H$  and
the number of quadrature points on $V$, $M$, tend to infinity.
%, and $\sigma^2$ the covariance of 
%the error term $ {\cal E}$ tends to zero. 
In practice, the error term $\delta$ is unknown, only $\tilde{\bu}$ is given:
  $F_{\bm, C}$ can be computed but not $F_{\bm, C}^\delta$.
	Note that $\delta $ tends to zero as $\dim H \ri \infty$.
	Therefore, we will from now on consider the computable  distribution function 
	\bean \label{real1}
\rho ( \tilde{\bu}|\bm,C) = 
(\f{3}{2 \pi e})^{\frac{3M}{2}} \prod_{j=1}^{M} (MC'(j))^{\f32} % \no \\
 (\det ( C^{-1}(\ca_{\bm}^{M })'\ca_{\bm}^{M } +   I))^{-\f12}
 F_{\bm, C}(\bg_{min})^{-\f{3M}{2}} ,
\eean
where this time  $\bg_{min}$ is the minimizer
of $F_{\bm, C}$ over $H$.  The resulting posterior distribution function is
\bean \label{real2}
\rho (\bm , C| \tilde{\bu}) \propto 
\rho (\tilde{\bu}|\bm, C) \rho_{pr} (\bm, C).
\eean
Choosing the priors of $\bm$ and $C$ to be independent 
and uniform in $B$ and in $[C_0, C_1]$ respectively, where $C_0$ and $C_1$
are positive constants, it follows from  (\ref{real1}, \ref{real2})
that 
\bean \label{real3}
\rho (\bm | \tilde{\bu}) \propto (\f{3}{2 \pi e})^{\frac{3M}{2}} \prod_{j=1}^{M} (MC'(j))^{\f32}
\int_{C_0}^{C_1}
(\det ( C^{-1}(\ca_{\bm}^{M })'\ca_{\bm}^{M } +   I))^{-\f12}
 F_{\bm, C}(\bg_{min})^{-\f{3M}{2}} dC .
\eean
	We will show that for all
$\bm$ in $B$ away from the true geometry parameter $\tilde{\bm}$,
the posterior marginal 
$\rho(\bm|\tilde{\bu} )$ 
%computed from (\ref{real1}-\ref{real2}) 
converges  to zero
as $M \ri \infty$ and $\dim H\ri \infty$. \\
%The following theorem explains in what sense the argument of the minimum of the functional 
%$\cf_{\bm,C}$ converges to the 	
	%slip solving the fault inverse problem, and how the argument of the minimum of $f_C$ converges to the geometry parameter
	%solving the fault inverse problem.
	As we shall see, to obtain convergence we are not just "adding" measurement points:
	this is why we now assume that there are $M_n$ measurement points, where
	$M_n$ is a strictly increasing  sequence  in $\NN$.
	The measurement points are  denoted in this section by
	$P_j^n, j=1, ..., M_n$ and in particular 
	the set $\{ P_j^n: j=1, ..., M_n \}$ is not necessarily included
	in the set $\{ P_j^{n+1}: j=1, ..., M_{n+1} \}$. 
	The coefficients $C'$ are now denoted by $C'(j,n)$, $j=1, ..., M_n$
	and are required to relate to a quadrature rule,  as specified below.
	Instead of just fixing a finite-dimensional subspace $H$ of $H^1_0(R)$,
	we have  to consider a sequence of finite-dimensional subspaces $H_n$. 
	%$$uch that 
	%$$H_N \subset H_{N+1}$, 
	%and $H_N $ may also depend on $\bm$.
	The operator defined in \eqref{AmpN def}
	is accordingly denoted by $\ca_{\bm}^{ n}$, as $M_n$ depends on $n$. 
	We make the following  assumptions on the sequence of finite subspaces $H_n$  and 
on the quadrature rule on $V$:
\begin{itemize}%{enumerate}%[label=H\arabic*. , wide=0.5em,  leftmargin=*]
 \item The quadrature with the weights $C'(j,n)$ and the nodes
 $P^n_j$, $j=1,..., M_n$ is of order 1, more precisely for all $\phi$
 in $C^1(\ov{V})$,
\bean \label{rule order}
|\sum_{j=1}^{M_n} C'(j,n) \phi (P^n_j) - \int_V \phi | =
O(\f{1}{M_n}) \sup_{V} |\nabla \phi| .
\eean  
	%where $D^2 \phi$ is the Hessian of $\phi$
	\item \bean C'(j,n) > 0
	\mbox{ and } C'(j,n) = O(\f{1}{M_n}), 
	\mbox{ uniformly in } j
	\label{Cjn assum}
	\eean
	\item
	\bean
	H_n \subset H_{n+1} \mbox{ and } \bigcup_{n=1}^\infty H_n
	\mbox{ is dense in } H^1_0(R)
	\label{HN assump}
	\eean
		\item The dimension of the space $H_n$ does not grow too fast relative to the number of measurement points $M_n$, more precisely, 
		\bean
		\dim H_n = O(M_n) \label{dim HN}
		\eean
%		\item $\dim H_p  \geq 3N$ and the $3N$ eigenvalues of the discrete operator
%		$\ca_{\bm}^{p, N}(\ca_{\bm}^{p, N })'  $ are equal to 
%		the first $3N$ eigenvalues of the symmetric and compact  operator
%		$\ca_{\bm }\ca_{\bm }'$.
\end{itemize}%{enumerate}
The existence of such spaces  $H_n$ 
and quadrature rules on $V$ 
is trivial if, for example, $R$ and $V$ are polygons.
$H_n$ can be constructed from spaces of finite elements. 
A possible choice for defining the points $P_j^n$ is
$$
  \{P_j^n: j=1, .., M_n \} = (\f{1}{n} \ZZ^2) \cap V.
$$
In that case $M_n$ is approximately equal to $n^2 |V|$, where $|V|$ is the surface area 
of $V$. A natural choice for the quadrature coefficients  $C'(j,n)$ is $C'(j,n) = \f{1}{n^2}$. 
In general, the set of 
points $\{P_j^n: j=1, .., M_n \}$ need not be laid on  such a regular grid.
%, see their
 %Theorem 4.5. 
%Their study but as mentioned in  section 7 of their 
%paper, it can be easily generalized to multidimensional integrals; in case of rectangular domains,
%the generalization is straightforward. 
%The estimate  \eqref{weight estimate} holds  if one starts 
%from the standard Gauss Legendre quadrature.

% note: even though we start from H_0^1(R) and not L^2(R), we still obtain 

%Gauss Legendre typically for basis of L^2
%to get basis for grad, H^1(R) take x_1 and x_2 derivative
%that's still in the basis because taking a derivative lowers the degree
%must be equal to zero on x_1=1, x_1=-1, so a subspace will do 

\subsection{The error functions $f^{disc}_C$ and the convergence of the arguments of their minima}

\begin{lem}\label{intermediate lemma}
Assume that $\ch_n$ converges weakly to $\ch$ in $H^1_0(R)$.
Fix $\bm$ in $B$.
Then $\ca_\bm \ch_n - \ca_\bm \ch $ converges uniformly to zero in $V$.
Let $\bm_n$ be a sequence in $B$ converging to $\bm$.
Then 
$\ca_{\bm_n} \ch_n - \ca_\bm \ch $ converges uniformly to zero in $V$.
\end{lem}
\textbf{Proof:}
According to  \eqref{Aabd},
\bean
&&|\ca_\bm \ch_n (x_1, x_2)- \ca_\bm \ch(x_1, x_2)|\no \\
&=&|\int_R  \bH_\bm(\bx, y_1, y_2)
 (\ch_n (y_1, y_2) - \ch (y_1, y_2))s d y_1 dy_2 | \no \\
&\leq& \sup_{ \bx \in V, (y_1,y_2) \in R } | \bH_\bm(\bx, y_1, y_2)s||R|^{\f12}
(\int_R  
 (\ch_n (y_1, y_2) - \ch (y_1, y_2))^2 d y_1 dy_2 )^{\f12},
\label{uni estimate}
\eean
and since $\ch_n$ converges strongly to $\ch$ in $L^2(V)$,
the first claim is proved.
To prove the second claim, it suffices to show that 
$\ca_{\bm_n}  \ch_n- \ca_\bm \ch_n$ converges uniformly to zero.
This is due to the estimate,
\bean
|\ca_{\bm_n} \ch_n (x_1, x_2)- \ca_\bm \ch_n(x_1, x_2)| \leq  \no \\
\sup_{ \bx \in V, (y_1,y_2) \in R } | \bH_{\bm_n}(\bx, y_1, y_2)s_{\bm_n} -
\bH_\bm(\bx, y_1, y_2)s_{\bm} ||R|^{\f12}
(\int_R  
 \ch_n (y_1, y_2) ^2 d y_1 dy_2 )^{\f12},
\label{uni estimate 2}
\eean
and the lemma is proved. $\Box$\\

Recall the definition \eqref{define bhdiscmc} of the function $f^{disc}_C$, which through $F_{\bm,C}$ also depends
 on $M_n$  and $H_n$.
It will be convenient in the proof of  the next lemma to distinguish the indices 
 of $M$ and $H$, so they will be denoted by  $M_n$ and $H_p$.
In addition,  to clarify that $ f^{disc}_C$ depends on $M_n$ and $H_p$, 
we will use the notation
$f^{disc}_{C,n,p}$
to  make the dependency explicit.
We also use the standard notation $B(\tilde{\bm} ,r)$ 
for the open ball in $\RR^3$ with center $\tilde{\bm} $ and radius $r$.
Roughly speaking, the following lemma expresses
that if the data is coming from the true value $\tilde{\bm}$,
then   for some $\eta>0$, if we  take $\bm$ in $B \setminus 
B(\tilde{\bm} ,\eta)$,  the error functional 
$f^{disc}_{C,n,p}(\bm)$ will be larger  compared to when we take $\bm$
close to $\tilde{\bm}$. 
The statement and the proof of the lemma are rather intricate because
we need this statement  to  be uniform as $n \ri \infty$, and $C \ri 0$.
\begin{lem}\label{second eta lemma}
Assume that 
	$\tilde{\bu} (P_j^n)= \ca_{\tilde{\bm}} \tilde{\ch}(P_j^n) + \tilde{\ce}(P_j^n)$,
	$j=1, .., M_n$ for
	some $\tilde{\bm}$ in $B$ and
	$\tilde{\ch} \neq 0$ in  $H_0^1(R)$. 
Fix $\eta>0$.
Then there exists   $C^*>0$,
such that for all $C_0$ in $(0, C^*)$, there exist
 $\eta'>0$, $\epsilon >0$, $n_0$, and $p_0$
such that,
\bean 
\sup_{ n>n_0, \dim_{H_p} >p_0}
\sup_{C_0 \leq C \leq 2 C_0,  \bm \in \ov{B(\tilde{\bm} , \eta' )}\cap B} 
f^{disc}_{C,n,p}(\bm) \no \\
<\inf_{ n>n_0, \dim_{H_p} >p_0}
\inf_{C_0\leq C \leq 2 C_0,  \bm \in B \setminus 
B(\tilde{\bm} ,\eta) } f^{disc}_{C,n,p}(\bm),  \label{imp estimate2}
\eean
if $\sum_{j=1}^{M_n}  C'(j,n) |\tilde{\ce}(P_j^n)|^2  < \epsilon$ for  $n>n_0$.
\end{lem}

\textbf{Proof:}
If we first  prove  \eqref{imp estimate2} in the case where
$\tilde{\ce}(P_j^n) =0$, for all $j=1, ..., M_n$,
let $ \sqrt{2 \epsilon}$ be equal
to
\bea
\inf_{ n>n_0, \dim_{H_p} >p_0}
\inf_{C_0\leq C \leq 2 C_0,  \bm \in B \setminus 
B(\tilde{\bm} ,\eta) } \sqrt{f^{disc}_{C,n,p}(\bm)} \\
- \sup_{ n>n_0, \dim_{H_p} >p_0}
\sup_{C_0 \leq C \leq 2 C_0,  \bm \in \ov{B(\tilde{\bm} , \eta' )}\cap B} 
\sqrt{f^{disc}_{C,n,p}(\bm)}.
\eea
Then, by the triangle inequality,  the strict  inequality \eqref{imp estimate2} 
will still hold if  $\sum_{j=1}^{M_n}  C'(j,n) |\tilde{\ce}(P_j^n)|^2  < \epsilon$ for  $n>n_0$.
Thus, without loss of generality, we can now assume that 
$\tilde{\ce}(P_j^n) =0$, for all $j=1, ..., M_n$ and  $n>n_0$.\\
 %since inequality \eqref{imp estimate2} 
%is strict
Arguing by contradiction, if \eqref{imp estimate2} 
 does not hold, then for  a sequence
$C_k$ converging to zero 
we have that
for all $\eta'>0$,  $n_0$, and $p_0$,
\bea
\sup_{ n>n_0, \dim_{H_p} >p_0}
\sup_{C_k \leq C \leq 2 C_k,  \bm \in \ov{B(\tilde{\bm} , \eta' )}\cap B} 
f^{disc}_{C,n,p}(\bm) \no \\
\geq \inf_{ n>n_0, \dim_{H_p} >p_0}
\inf_{C_k \leq C \leq  2 C_k,  \bm \in B \setminus 
B(\tilde{\bm} ,\eta) } f^{disc}_{C,n,p}(\bm), % \label{cont}
\eea
and as $f^{disc}_{C,n,p}$ is increasing in $C$, %after possibly extracting a subsequence
\bean
\sup_{ n>n_0, \dim_{H_p} >p_0}
\sup_{ \bm \in \ov{B(\tilde{\bm} , \eta' )}\cap B} 
f^{disc}_{2C_k,n,p}(\bm) %\no \\
\geq \inf_{ n>n_0, \dim_{H_p} >p_0}
\inf_{ \bm \in B \setminus 
B(\tilde{\bm} ,\eta) } f^{disc}_{C_k,n,p}(\bm).  \label{cont1}
\eean
We can also assume that $\eta'=\eta_k =O(C_k^{\f12})$.
We set $n_0=a_k, p_0=b_k$
where  the sequences $a_k$ and $b_k$ are strictly increasing in $\NN$.
%that diverge to infinity and a sequence $u_n$ in $\RR$ converging to zero, 
%\bean
%\sup_{ \bm \in \ov{B(\tilde{\bm} , \eta_n)}\cap B} 
%f^{disc}_{2C_n,a_n}(\bm) 
%\geq u_n+
%\inf_{ n>a_n, \dim_{H_p} >b_n}
%\inf_{ \bm \in B \setminus 
%B(\tilde{\bm} ,\eta) } f^{disc}_{C_n,n}(\bm), \label{cont1ver2}
%\eean
%where $f^{disc}_{C,a_n}(\bm)$ is the minimum of the $F_{\bm, C}$
%that  uses the points $P_j^{M_{a_n}}$, 
%$j=1, ..., M_{a_n}$ and the finite-dimensional space 
%$H_{b_n}$.
Let $\Pi_{b_k}$ be the orthogonal projection on $H_{b_k}$.
As $\Pi_{b_k} \tilde{\ch} $ converges to $\tilde{\ch}$
due to \eqref{HN assump},
$\ca_{\tilde{\bm}} \Pi_{{b_k}} \tilde{\ch} $
converges uniformly to
$\ca_{\tilde{\bm}} \tilde{\ch}$  in $V$
by Lemma \ref{intermediate lemma}.
After possibly extracting a   subsequence, we can assume that 
\bean
 \sup_{(x_1, x_2) \in V}
|(\ca_{\tilde{\bm}} \Pi_{b_k} \tilde{\ch} )(x_1, x_2) -
(\ca_{\tilde{\bm}} \tilde{\ch} )(x_1, x_2)  | 
= O(C_k^{\f12}) .
\label{uni on V 0}
\eean
%and that $u_n = O(C_n) $.
%Since 
%$f^{disc}_{C,n}$ is increasing in $C$, %after possibly extracting a subsequence
%\bean 
%\sup_{ \bm \in \ov{B(\tilde{\bm} , \eta_n)}\cap B} 
%f^{disc}_{2C_n,a_n}(\bm) 
%\geq 
%\inf_{\bm \in B \setminus 
%B(\tilde{\bm} ,\eta) } f^{disc}_{C_n,a_n}(\bm), %\label{cont2}
 %\label{cont2}
%\eean
%while (\ref{uni on V 0})   still holds. 
By continuity 
$\ds \sup_{  \bm \in \ov{B(\tilde{\bm} , \eta_k )}\cap B} 
f^{disc}_{2C_k, n,p}(\bm) $
 is achieved at some $\bm_k$ in 
$\ov{B(\tilde{\bm} , \eta_k)}$. 
%Note that by construction 
%$|\bm_n - \tilde{\bm} |=O(C_n^{\f12})$. 
It is also bounded above by
\bean \label{upper bound}
\sum_{j=1}^{M_{n}} C'(j,n) | 
(\ca_{\bm_k} \Pi_{p}\tilde{\ch}  - \ca_{\tilde{\bm}} \tilde{\ch} |^2  (P_j^{n})
+ 2C_k \int_R |\nabla \Pi_{p}\tilde{\ch}  |^2.
\eean
Given that $|\bm_k - \tilde{\bm} |=O(C_k^{\f12})$, 
we can apply 
(\ref{Cjn assum}, \ref{HN assump}, \ref{uni estimate 2}, \ref{uni on V 0})
%,\ref{uni on V 2},\ref{uni on V}) 
to find that for all $n>a_k$, $p>b_k$,
\bea
\sum_{j=1}^{M_{n}} C'(j,n) | 
(\ca_{\bm_k} \Pi_{p}\tilde{\ch}  - \ca_{\tilde{\bm}} \tilde{\ch} |^2  (P_j^{n})
= O(C_k).
\eea
%for all $N>N_n$. 
%Recalling $\epsilon_n= O(C_n)$,
It follows that 
%the expression in 
%\eqref{upper bound} is $O(C_n)$ and so is 
the left hand side of \eqref{cont1} is $O(C_k)$, if $n_0=a_k$, $p_0=c_k$.
From the right hand side of  \eqref{cont1},  we now claim that 
there exist   a sequence
$\bm_k'$ in $B\setminus B(\tilde{\bm}, \eta)$
and $a_k'$ and $b_k'$ two strictly increasing sequences in $\NN$,
such that $
f^{disc}_{C_k, a_k', b_k'}(\bm_k') = O(C_k) 
$.  
$f^{disc}_{C_k, a_k', b_k'}(\bm_k') $ is equal to
\bean \label{the inf}
\sum_{j=1}^{M_{a_k'}} C'(j,a_k') | 
(\ca_{\bm_k'} \bh_k  - \ca_{\tilde{\bm}} \tilde{\ch} |^2  (P_j^{a_k'})
+ C_k\int_R |\nabla \bh_k  |^2,
\eean
for some $\bh_k$  in $H_{b_k'}$.
Since the quantity in \eqref{the inf} is $O(C_k)$,  
$\bh_k $ is  bounded in $H^1_0(R)$, so after possibly  extracting a subsequence,
we may assume that it is weakly convergent to some 
$\bh^*$ in $H^1_0(R)$.
Likewise, 
we may assume that $\bm'_k$ converges to some
$\bm^*$ in $B\setminus B(\tilde{\bm}, \eta)$.
From  Lemma \ref{intermediate lemma}  we can claim that,
\bea
\sum_{j=1}^{M_{a_k'}} C'(j,a_k') | 
\ca_{\bm_k'} \bh_k  - \ca_{\bm^*} \bh^*|^2  (P_j^{a_k'})
\eea
converges to zero. Recalling 
\eqref{the inf},
this shows that
\bea
\sum_{j=1}^{M_{a_k'}} C'(j,a_k') | 
\ca_{\bm^*} \bh^* - \ca_{\tilde{\bm}} \tilde{\ch} |^2  (P_j^{a_k'})
\eea
converges to zero.
But this term also converges to $\int_V | \ca_{\bm^*}\bh^*
- \ca_{\tilde{\bm}} \tilde{\ch}|^2$, contradicting the injectivity 
statement in Proposition 
\ref{injectivity of Am} 
as $|\bm^* - \tilde{\bm}|\geq \eta >0$ and $\tilde{\ch}  \neq 0$.
$\Box$

% tasks:
% 1. rephrase this theorem in the m context
% 2. Are C_0 and C_1 necessary?

\subsection{A uniformly bounded determinant}
\begin{lem} \label{detlem}
Fix two positive constants $C_0$ and $C_1$ such that $C_0<C_1$.
If \eqref{rule order} and \eqref{dim HN} hold then
the determinant $\det (C^{-1} (\ca_{\bm}^{n })' \ca_{\bm}^{n} 
 + I)$ is bounded below by 1 and above by a constant, 
for all $n$ in $\NN$, $\bm$ in $B$, and $C$ in $[C_0, C_1]$.
\end{lem}
\textbf{Proof:}
Using an orthonormal basis of $H_n$ made of eigenvectors 
$\phi_1, ..., \phi_{q_n}$
of 
$(\ca_{\bm}^{ n })' \ca_{\bm}^{ n} $,  we have  the estimate
\bea
\det (C^{-1} (\ca_{\bm}^{n })' \ca_{\bm}^{n} 
 + I) \leq
\exp ( C^{-1} \sum_{j=1}^{q_n} \lambda_j) \leq
 \exp ( C_0^{-1} \sum_{j=1}^{q_n} \lambda_j) ,
\eea
where $(\ca_{\bm}^{ n })' \ca_{\bm}^{n} \phi_j = \lambda_j \phi_j$ and
$q_n$ is the dimension of $H_n$.
%We use $\Pi_p$
%Let $Q_N$ be the orthogonal  defined on $L^2(V)$ with range
%the finite-dimensional space
%$\ca_\bm H_p$.
%Let $Q_N^{disc}$ be the linear operator from $C^2(V)^3$ to 
%$\RR^{3 M_N}$
%\bea
%Q_N^{disc} (\phi) =
%(C'(1,N) \phi(P_1^N), ..., C'(M_N,N) \phi(P_{M_N}^N))
%\eea 
%Let $\{ \Psi_j \in L^2(V): j=1, ..., M_N \}$ be an orthonormal set and define
%the linear operator %orthogonal projection
%\bea
%D_N: C^2(\ov{V})^3 \ri L^2(V)^3,\\
%D_N(\varphi) = \sum_{j=1}^{M_{N_n}} \sum_{k=1}^{3}
%C'(j,N)^{\f{1}{2}}( \varphi(P_j^N) \cdot \bev_k) \bev_k \Psi_j.
%\eea 
%We note that for all $\bg, \bh$ in $H_p$,
%\bean
%<\ca_\bm^N \bg, \ca_\bm^N \bh> = <D_N \ca_\bm\bg, D_N \ca_\bm \bh>, 
%\label{inner products}
%\eean
%where the first inner product is the natural dot product in $\RR^{3M_N}$
%and the second inner product the natural inner product in $L^2(V)^3$.\\
%We now estimate the difference 
%$<D_N \ca_\bm\bg, D_N \ca_\bm \bh> - <\ca_\bm\bg, \ca_\bm \bh> $ 
%where $\bg$ and $\bh$ are in $H_p$
%and both inner products are those of $L^2(V)^3$.
%norm of $D_N \ca_\bm - \ca_\bm$.
Let $\bg$ be in $H^1_0(R)$. %
 %such that $\| \bg \| =1$. 
 Given the definition of $\ca_\bm$ \eqref{Aabd}
and that $\bH_\bm(\bx, y_1, y_2)$ and
$\p_{x_i} \bH_\bm(\bx, y_1, y_2)$, $i=1,2$, %
%$\p x_i\p x_j \bH_\bm(\bx, y_1, y_2)$, 
are continuous in $\bx, (y_1, y_2), \bm$ in $\ov{V}\times\ov{R} \times B$,
it follows that 
$$ \sup_{\bx \in \ov{V}} |\ca_\bm  \bg| +
\sup_{\bx \in \ov{V}} |\nabla \ca_\bm  \bg| 
%\sup_{\bx \in \ov{V}} |D^2 \ca_\bm  \bg|
 = O(\|  \bg \|), $$
uniformly for $\bm$ in $B$.
Therefore, from \eqref{rule order}, 
 if $\bg, \bh$ are in $H^1_0(R)$ such that $\| \bg \|= \| \bh \| =1$,
\bean
 &&\sum_{j=1}^{M_{n}} 
C'(j,n) \ca_\bm  \bg(P_j^n) \cdot \ca_\bm  \bh(P_j^n) 
- \int_V \ca_\bm  \bg  \cdot \ca_\bm  \bh \no \\
&=& O(\f{1}{M_n}) \sup_{V} |\nabla (\ca_\bm  \bg \cdot \ca_\bm  \bh)|\no \\
&=&  O(\f{1}{M_n}),  \label{disc est}
\eean
unformly for $\bm$ in $B$.
%Thus,
%\bean
%<D_N \ca_\bm\bg, D_N \ca_\bm \bh> - 
%<\ca_\bm\bg, \ca_\bm \bh>=  O(\f{\|  \bg \| |\bh \|}{M_N^2}),
%\label{OghM}
%\eean
%for all $\bg, \bh$ in $H^1_0(R)$. 
But, 
\bea
\sum_{k=1}^{q_n} \lambda_k
=\sum_{k=1}^{q_n} \sum_{j=1}^{M_{n}}
C'(j,n) \ca_\bm  \phi_k(P_j^n) \cdot \ca_\bm  \phi_k(P_j^n).
\eea
Recalling \eqref{dim HN}, $ q_n = O(M_n)$, so by \eqref{disc est}, 
\bea
\sum_{k=1}^{q_n} \lambda_k
=\sum_{k=1}^{q_n}  \int_V \ca_\bm \phi_k   \cdot \ca_\bm  \phi_k + O(1).
\eea 
%is also the trace of 
%$(\ca_{\bm}^{N })' \ca_{\bm}^{N} $ and due to \eqref{inner products}
%it is also the trace of the finite range operator 
%$(D_N \ca)^* D_N \ca$.
Since $\ov{R}$ and $\ov{V}$ are compact 
and $\bH_\bm$ is smooth on $\ov{V} \times \ov{R}$,
given the definition \eqref{Aabd} of $\ca_\bm$, $\ca_\bm^* \ca_\bm$ is a trace class operator
\cite{gohberg1996traces}.
%Now let $\mu_i$
The term $ \sum_{k=1}^{q_n} 
<\ca_\bm \phi_k , \ca_\bm  \phi_k >$ is bounded above by
$\|  \ca_m^* \ca_m \|_1 =\|  \ca_m \ca_m^* \|_1$ 
where $\| \, \|_1$ is the trace class norm (\cite{gohberg1996traces}, Theorem 3.5, Chapter
IV).
Thus there only remains to prove $\|  \ca_m \ca_m^* \|_1$
is uniformly bounded for $\bm$ in $B$.\\
Recalling the definition of $\ca_\bm $  \eqref{Aabd}, 
we find that $\ca_\bm \ca_\bm^* $ can be given in integral form,
\bea
\ca_{\bm} \ca_\bm^* &:& L^2(V)^3 \ri L^2(V)^3 \no \\
&&\bu \ri \int_V K_\bm(\bx, \bx')
 \bu (\bx')  d \bx' , %.  \label{Aabd}
\eea
where 
\bea
K_\bm(\bx, \bx') =\int_R  \bH_\bm(\bx, y_1, y_2) \bH_\bm'(\bx', y_1, y_2)
  s^2 d y_1 dy_2 .
\eea
Given that $K_\bm$ is continuous in $(\bx,\bx')$, that $V$ is bounded
 and $\ca_\bm \ca_\bm^* $ is  trace class, we have the following explicit 
formula \cite{gohberg1996traces}, Theorem 8.1, Chapter IV, 
\bea
\mbox{tr }(\ca_\bm \ca_\bm^*) = \mbox{tr }
\int_V  K_\bm(\bx, \bx) d \bx = \mbox{tr }
\int_V  \int_R \bH_\bm(\bx, y_1, y_2) \bH_\bm'(\bx, y_1, y_2) 
  s^2 d y_1 dy_2 d \bx.
\eea
As $\bH_m(\bx, y_1, y_2)
  s$ is uniformly bounded for $\bx$ in $V$, $\by$ in $R$, and $\bm$ in $B$,
	the result follows. $\Box$\\

\subsection{Convergence  of the posterior of $\bm$
as $M_n$ and $\dim H_p$ tend to infinity}

\begin{thm} \label{main cv theorem}
Assume that the 
data $\tilde{\bu}$ is  given
by (\ref{noise model}) with the true values $\bm = \tilde{\bm}$,
$\cg = \tilde{\ch}$, $\bg = \Pi_n \tilde{\ch}$, 
and is such that the 
random variables
$\ce(P_j^n) \cdot \bev_k$ where $ j=1, .., M_n$, $k=1,2,3$, $n \in \NN$,
%and $M_N$ increases to infinity,
are independent and identically distributed with zero mean and finite covariance $\sigma^2 $.
%Assume that the random variables
%$\ce(P_j^N) \cdot \bev_k$ where $ j=1, .., M_N$, $k=1,2,3$, $N \in \NN$,
%and $M_N$ increases to infinity,
%are independent and identically distributed with finite covariance $\sigma^2 $.
% normally distributed with  zero mean and covariance $\sigma^2 $.
Suppose conditions (\ref{rule order}-\ref{dim HN}) hold 
and that $\bm$ follows the distribution \eqref{real3}.
	Fix $\eta >0$.
Then there for all  $C_0>0$ small enough there is a positive $\sigma_0 $ such that if $\sigma < \sigma_0$,  
the probability of $\{  \bm \in B: |\bm - \tilde{\bm}| > \eta\} $
	converges  to zero as $n$  tends to infinity. 
%	$\bm $ in $B \setminus 
%B (\tilde{\bm}, \eta)$ as $N$ and $p$ tend to infinity.
% and $\sigma$
%tends to zero.
\end{thm}
\textbf{Proof:}
%Must use $C^*$.\\
We note that due to \eqref{Cjn assum} there is a positive constant $\alpha$
such that,
\bea
\sum_{j=1}^{M_n}  C'(j,n) |\ce (P_j^n)|^2  \leq 
\f{\alpha}{M_n}
\sum_{k=1}^3 \sum_{j=1}^{M_n}  (\ce (P_j^n)\cdot \bev_k)^2 . 
\eea
%By assumption, the $|\tilde{\ce}(P_j^N)|^2 $ are independent with expectation
%$3 \sigma^2$.
Using the law of large numbers (\cite{girardin2018applied}, Theorem 1.93),
$\f{1}{M_n} \sum_{k=1}^3 \sum_{j=1}^{M_n}   (\ce (P_j^n)\cdot \bev_k)^2 $
%follows a Gamma distribution with parameters $\alpha = \beta = 3 M_N/2$. 
%Therefore as $M_N \ri \infty$,  $\f{1}{3 M_N} \sum_{k=1}^3 \sum_{j=1}^{M_N}   \f{1}{\sigma^2}(\tilde{\ce}(P_j^N)\cdot \bev_k)^2 $ converges to 1 in probability. 
converges almost surely to $\sigma^2$ as $n \ri \infty$.
Following Lemma \ref{second eta lemma}, let $\epsilon >0$ such that inequality
\eqref{imp estimate2} is satisfied and let $\sigma_0 $ be  a fixed positive number less than
 $\f{\epsilon}{4 \alpha} $. 
Then $\sum_{j=1}^{M_n}  C'(j,n) |\ce (P_j^n)|^2 < \f{\epsilon}{4} $
almost surely  as $n\ri \infty$.  
Without loss of generality, we may assume that $n$ 
is large enough so that
$\sum_{j=1}^{M_n} C'(j,n) | \delta (P_j^n)|^2 < \f{\epsilon}{4} $,
since $\delta = \ca_\bm(I-\Pi_n) \tilde{\ch}$.
Then almost surely, as $n\ri \infty$,
$$\sum_{j=1}^{M_n}  C'(j,n) |(\ce  + \delta) (P_j^n)|^2 < \epsilon .$$
Now let $C^*, \eta', n_0, p_0  $ be as in the statement of Lemma \ref{second eta lemma} and assume that $n$ is such that $n>n_0$ and $\dim H_n > p_0$ while
$C_0<C^*$. Let $\ov{C}$ be such that 
\bea
  [C_0, 2 C_0] \cap [C_0, C_1] = [C_0, \ov{C}].
\eea
According to (\ref{real3}), 
\bea
\rho (\bm | \tilde{\bu})= {\cal I} {\cal C}_n \int_{C_0}^{C_1} 
 (\det ( C^{-1}(\ca_{\bm}^{n })'\ca_{\bm}^{ n } +   I))^{-\f12}
 F_{\bm, C}(\bg_{min})^{-\f{3M_n}{2}}
d C ,
\eea
where
\bea
{\cal C}_n = (\f{3}{2 \pi e})^{\frac{3M_n}{2}} \prod_{j=1}^{M_n} (M_nC'(j,n))^{\f32},
\eea
and
\bea
{\cal I}^{-1} =  {\cal C}_n \int_{B} \int_{C_0}^{C_1} 
 (\det ( C^{-1}(\ca_{\bm}^{n })'\ca_{\bm}^{ n } +   I))^{-\f12}
 F_{\bm, C}(\bg_{min})^{-\f{3M_n}{2}}
d C d \bm.
\eea
Set 
$\gamma_1$ to be the left hand side of \eqref{imp estimate2} and
$\gamma_2$ to be the right hand side of \eqref{imp estimate2}.
According to Lemmas \ref{second eta lemma} and \ref{detlem},
for $n>n_0$ and $\dim H_n > p_0$,
\bea
{\cal I}^{-1} &\geq &
 {\cal C}_n O(  \int_{B} \int_{C_0}^{C_1}   F_{\bm, C}(\bg_{min})^{-\f{3M_n}{2}}
d C d \bm ) \\
&\geq &{\cal C}_n O(   \int_{ B(\tilde{\bm} , \eta' )\cap B}   \int_{C_0}^{\ov{C}}
 F_{\bm, C}(\bg_{min})^{-\f{3M_n}{2}}
d C d \bm )\\
&\geq &{\cal C}_n O(     \int_{  B(\tilde{\bm} , \eta' )\cap B } \int_{C_0}^{\ov{C}}
\gamma_1^{-\f{3M_n}{2}}
d C d \bm ),
\eea
thus
\bea
{\cal I} \leq {\cal C}_n^{-1} O(\gamma_1^{\f{3M_n}{2}}).
\eea
Now assume that $\bm $ is in $B \setminus 
B(\tilde{\bm} ,\eta)$. 
If $C $ is in $[C_0, \ov{C}]$, then $F_{\bm, C}(\bg_{min}) \geq \gamma_2$.
If $C$ is in $[\ov{C}, C_1]$ then as $F_{\bm, C}$ is increasing in $C$, we still have
that $F_{\bm, C}(\bg_{min}) \geq \gamma_2$.
It follows that
\bea 
\rho (\bm | \tilde{\bu}) \leq {\cal I} {\cal C}_n O(\gamma_2^{-\f{3M_n}{2}})
\leq O((\f{\gamma_1}{\gamma_2})^{\f{3M_n}{2}}).
\eea
%uniformly for all $C$ in $[C_0, C_1]$.
 %and $n> n_0$.  
%This is also true for $C$ in $[C^*, C_1]$ since $F_{\bm, C}(\bg_{min})=
%f^{disc}_C(\bm)$ is increasing in $C$. 
%
% large enough in $\NN$. 
%As  the volume of  $\{  \bm \in B: |\bm - \tilde{\bm}| > \eta\} $ is finite,
We have found that the probability of $\bm \in \{  \bm \in B: |\bm - \tilde{\bm}| > \eta\} $ is $ O((\f{\gamma_1}{\gamma_2})^{\f{3M_n}{2}})$,
if $\ce$ satisfies
$\sum_{j=1}^{M_n}  C'(j,n) |\ce (P_j^n)|^2  < \f{\epsilon}{4}$ for $n>n_0$. 
As by Lemma \ref{second eta lemma}, 
$0< \f{\gamma_1}{\gamma_2} < 1$,  the proof is complete. $\Box$ \\
\textbf{Remark:}
Let $\Omega$ be the underlying probability space for the measurements at the points 
$P_j^n$. Theorem \ref{main cv theorem} assumes that the  random variables
$\ce (P_j^n) \cdot \bev_k$, $ j=1, .., M_n$, $k=1,2,3$, $n \in \NN$
where $M_n$ is increasing and tends to infinity, are independent and identically
distributed. The distribution \eqref{full law} was constructed under the assumption 
that $\ce (P_j^n) $ is Gaussian, however, once this  distribution is set,
this is no longer necessary for the convergence result of 
Theorem \ref{main cv theorem} to hold.
In practice, Theorem \ref{main cv theorem} expresses that 
if measurements $\ce (P_j^n)(\omega) \cdot \bev_k$, $ j=1, .., M_n$, $k=1,2,3$, $n \in \NN$
are available, almost surely for $\omega $ in $\Omega$,
if $\sigma < \sigma_0$, 
the probability of $\{  \bm \in B: |\bm - \tilde{\bm}| > \eta\} $
is $O((\f{\gamma_1}{\gamma_2})^{\f{3M_n}{2}})$. %, if $p > p_0$.
%$\omega, \Omega$ sequence of measurements $\tilde{\ce}(P_j^N) (\omega)$, $\omega \in \Omega$,
%$j=1, .., M_{N_n}$, $N_n$ increasing to infinity

\section{Extension to more general inverse problems}\label{extension}
The reconstruction algorithm developed in this paper 
was initially designed
to solve a specific fault inverse problem arising in geophysics.  
%and the proof of the convergence of its
%solution to the solution of the underlying continuous problem were initially designed
%to solve a specific fault inverse problem arising in geophysics  
However, it can be used in other inverse problems where both linear and nonlinear unknowns have to be recovered. This algorithm is especially well-suited to problems where the linear part has to be regularized, but to what extent is unknown, and this lack of knowledge is made more challenging 
by the fact that the underlying linear operator depends on the linear unknown. 
This algorithm is also  well-suited to problems where noise levels are considerable and
as a result
 the maximum 
likelihood solution may become  meaningless. 
The expected value of the solution and its 
(rather large) covariance are then more meaningful. 
Section 4.3 of \cite{volkov2020stochastic} shows how 
this method outperforms
other regularization parameter selection methods
such as GCV or ML (used globally over the range of the nonlinear parameter)  on such problems.
The numerical simulations are considerably more involved in \cite{volkov2020stochastic} since the nonlinear parameter there is in $\RR^6$ and thus the gains made by applying the method advocated
in this paper are even more dramatic.\\
%Introduce lots of general notations 
We now write down a more general inverse problem that can efficiently be
solved following the algorithm analyzed in this paper.
For ease of notations, we only consider the case of scalar functions:
the case of vector fields can easily be inferred from there. 

\subsection{Inverse problem formulation}
Let $R$ and $V$ be two open bounded subsets of $\RR^d$. 
Let the operator
\bean \label{generalAm}
\ca_{\bm} &:& L^2 (R) \ri L^2(V) 
\eean
be a compact linear operator that depends continuously on $\bm$, a parameter
in $B \subset \RR^q$. We assume that $B$ is compact.\\\\
Uniqueness assumption: \\
\textit{
For any $\bm_1$ and $\bm_2$ in $B$ and any $\cg_1$ and $\cg_2$ in
$H^1_0 (R)$, if  $\ca_{\bm_1} \cg_1 = \ca_{\bm_2} \cg_2 $ in $L^2(V)$,
then $\bm_1=\bm_2$ and $\cg_1=\cg_2$.}

\subsection{Related discrete inverse problem}
Let $P_j$,  $j=1, ..., M$ be points in $V$, referred to as measurement points.
The discrete inverse problem consists of finding $\bm$
from noisy data at the points $P_j$.
More precisely, let $H$ be a finite-dimensional subspace of $H^1_0(R)$.
%$p$ is an index $H_p$  may 
We assume that the noisy data is given in the form
\bean% \label{noise model}
 %m = \argmin_{ m \in B} \argmin_{\bg \in H^1_0(R)} F^{disc}_{m,C} (\bg)
(\tilde{\bu} (P_1), ..., \tilde{\bu} (P_{M}) ) =  
(\ca_{\tilde{\bm}} \tilde{\ch} (P_{1}), ..., \ca_{\tilde{\bm}}  \tilde{\ch}(P_{M}) )
+ {\cal E}, \label{data def}
\eean
where ${\tilde{\bm}}\in B$,  $\tilde{\ch} \in H^1_0 (R)$.
Note that earlier in this paper $(\ca_{\tilde{\bm}} \tilde{\ch} (P_{1}), ..., \ca_{\tilde{\bm}}  \tilde{\ch}(P_{M}) )$  had to be split in two terms, 
one corresponding to the orthogonal projection $\Pi \tilde{\ch} $  
on $H$ and the other term, $\delta$, related to the residual
$(I-\Pi) \tilde{\ch} $. This splitting was instrumental to prove the convergence result,
but is no longer necessary to formulate the posterior distribution of $\bm$
or to state the convergence result.
Here, ${\cal E}$ is a 
random variable with mean zero and covariance $\sigma^2 I$.
%It is in fact not necessary to assume that ${\cal E}$ 
%is normal, neither for constructing
%the posterior distribution of $\bm$ nor for proving 
%that it converges to $\tilde{\bm}$.
%The term $\delta $ is an error term due to reduction to the 
%finite dimensional space $H_p$, it is in the form
%\bea
%\label{delta def}
%(\ca_\bm (I- \Pi_p)(\cg) (P_1^N), ..., \ca_\bm (I- \Pi_p)(\cg) (P_{M_N}^N) ),
%\eea
%where $\Pi_p$ is the orthogonal projection from 
%$H^1_0(R) $ to $H_p$ and $\cg$ is in $H^1_0(R) $.
The vector $(\tilde{\bu} (P_1), ..., \tilde{\bu} (P_{M}) ) $
will be denoted by $\tilde{\bu}$.
\\\\
Discrete inverse problem statement:\\
\textit{Assume that the sets $R, V, B$ are known. Assume that the mapping
$\bm \ri \ca_{\bm}$ is known.
Recover $\bm$ from $\tilde{\bu} $: compute its expected value and its covariance, and if possible, 
its marginal probability distribution functions.
} \\\\
Let $C'(j)$ be positive weights.
%for a quadrature rule on $V$ corresponding to the points 
%$P_j^N$.  
For a fixed $C>0$, a fixed $\bm$ in $B$, 
and $\bg$ in $H$,  define the error functional
\bean
F_{\bm,C} (\bg) = \sum_{j=1}^{M} C'(j) | 
(\ca_{\bm} \bg- \tilde{\bu} ) (P_j)|^2 
+ C \int_R |\nabla \bg|^2 .
\label{Ffunc def}
\eean
To ease notations, it is convenient to introduce the finite-dimensional operator 
$\ca_{\bm}^{M}$ which is built as  a slight correction of $\ca_\bm $ by the quadrature
coefficients $C'(j)$:
\bean
\ca_{\bm}^{M} \bg =
(C'(1)^{\f12}  \ca_\bm (\bg)(P_1), ..., C'(M)^{\f12}  \ca_\bm (\bg)(P_{M})).
\label{AmpN def2}
\eean
Define the probability
distribution
\bean \label{rho def}
\rho (\bm , C| \tilde{\bu}) \propto 
(\det ( C^{-1}(\ca_{\bm}^{M })'\ca_{\bm}^{M } +   I))^{-\f12}
 F_{\bm, C}(\bg_{min})^{-\f{M}{2}}
 \rho_{pr} (\bm, C),
\eean
where $\propto $ means "proportional to", $F_{\bm, C}(\bg_{min})$ is 
the minimum of  $F_{\bm,C}$ over $H$,
and $ \rho_{pr} (\bm, C)$ is the prior distribution of $(\bm, C)$
(note that the exponent of $F_{\bm, C}(\bg_{min})$ is different 
from the one in formula \eqref{real1}: this is because in the case of the elasticity 
equation we have $M$ measurements of  three-dimensional vectors, while here
we assumed that we have $M$ scalar measurements).
If we choose the priors of $\bm$ and $C$ to be independent 
and uniform in $B$ and in $[C_0, C_1]$ respectively, where $C_0$ and $C_1$
are positive constants, it follows from  \eqref{rho def}
that 
\bean \label{rho def 2}
\rho (\bm | \tilde{\bu}) \propto \int_{C_0}^{C_1}
(\det ( C^{-1}(\ca_{\bm}^{M })'\ca_{\bm}^{M } +   I))^{-\f12}
 F_{\bm, C}(\bg_{min})^{-\f{M}{2}} dC .
\eean
Note that another 
typical choice for the prior distribution of $ C$ is to take it  
such that $\log_{10} C$
 is uniformly distributed on an interval in $\RR$.\\
Regarding the computation of $\rho (\bm , C| \tilde{\bu})$, as
$\bm \in B \subset \RR^q$, if $q=1$
 then the  trapezoidal rule is adequate. However, if $q \geq 3$, 
Markov Chain Monte Carlo techniques are called for. If a parallel computing platform
is available, we recommend using the parallel sampling algorithm for inverse problems
combining linear and nonlinear unknowns described in  \cite{volkov2020stochastic}, or any algorithm
along these lines.

\subsection{Convergence of the solution to the discrete inverse problem to  the unique solution of the underlying continuous problem}
The convergence Theorem \ref{main cv theorem} can easily be extended 
to the more general inverse problem considered in this section. 
%If $\tilde{\bm}$ is the true value 
%We can prove that the probability density  
For the arguments used in the proof of Lemmas 
\ref{intermediate lemma}, \ref{second eta lemma}, and 
\ref{detlem} to be valid, 
the operator \eqref{generalAm}
has to be associated with an integration kernel 
$\bH_\bm$ such that 
\bean \label{Am1}
 \ca_\bm \bg (\bx) =  \int_R \bH_\bm(\bx, \by)
 \bg (\by)  d \by, 
\eean
 for all $\bg$ in $L^2(R)$ and where $\bH_\bm$ is such that
\bean \label{Am2}
(\bm,\bx,\by ) \ri \bH_\bm(\bx, \by) \mbox{ is continuous  in }
B \times \ov{V} \times \ov{R}, 
\eean
and
\bean \label{Am3}
(\bm,\bx,\by ) \ri \p_{x_i}\bH_\bm(\bx, \by) \mbox{ is continuous  in }
B \times \ov{V} \times \ov{R}, \, i =1, ..., d. 
\eean
We are able to prove a convergence result if 
 the coefficients $C'$ come from a quadrature rule satisfying
(\ref{rule order}, \ref{Cjn assum})
 and if the finite-dimensional space $H=H_n$ 
satisfies (\ref{HN assump}, \ref{dim HN}).
Note that if the geometries of $R$ and $V$ are simple enough (for examples
if they are polygons), these conditions can be easily met.

\begin{thm} \label{general cv theorem}
Let the 
data $\tilde{\bu}$ be  given
by \eqref{data def} with the true values $\tilde{\bm}$,
$\tilde{\ch}$,  where $\ca_\bm$ satisfies (\ref{Am1}-\ref{Am3}),
and ${\cal E}$  a random vector such that its coordinates
are independent and identically distributed with zero mean and finite covariance $\sigma^2 $.
%for all $ j=1, .., M_N$ and $N \in \NN$.
%Suppose that $\bm$ follows the distribution 
%$\rho (\bm | \tilde{\bu})$ derived as the marginal of
  %(\ref{rho def}) through (\ref{Ffunc def}, \ref{AmpN def2}).
	%Assume also that the priors of $\bm$ and $C$ are independent, and that 
%the prior of $C$ is the uniform distribution  in $[C_0, C_1]$, where $C_0>0$,
%while the prior of $\bm$ is the uniform distribution  in $B$.
Assume that the coefficients $C'$ come from a quadrature rule satisfying
(\ref{rule order}, \ref{Cjn assum})
 and consider a sequence of 
 finite-dimensional spaces  $H$
satisfying (\ref{HN assump}, \ref{dim HN}).
Fix a distance $\eta>0$.
	Then there for all  $C_0>0$ small enough there is a positive $\sigma_0 $ such that if $\sigma < \sigma_0$,  
the probability of $\{  \bm \in B: |\bm - \tilde{\bm}| > \eta\} $
under the distribution
given by \eqref{rho def 2}
 converges  to zero as $M$ and $\dim H$ tend to infinity.
	\end{thm}

\section{Numerical simulations}\label{Numerical simulations}
% codes are in
% C:\Users\darko\Desktop\RESEARCH\Projects for 2020\calc for third summer paper
% there are six directories corresponding to the six cases: 12,12,25,25,50,50
% in each directory
% data is prepared in build_data_part1.m
% inverse problem is solved in par_inversion_august2020.m
% finally for producing graphs for solutions,
% post_process.m
% and 
% show_recon_slip.m
% see also directory 'plot pdfs together'
% for plotting pdfs together
\subsection{Construction of data}
\begin{figure}[htbp]
    \centering
      \includegraphics[scale=.4]{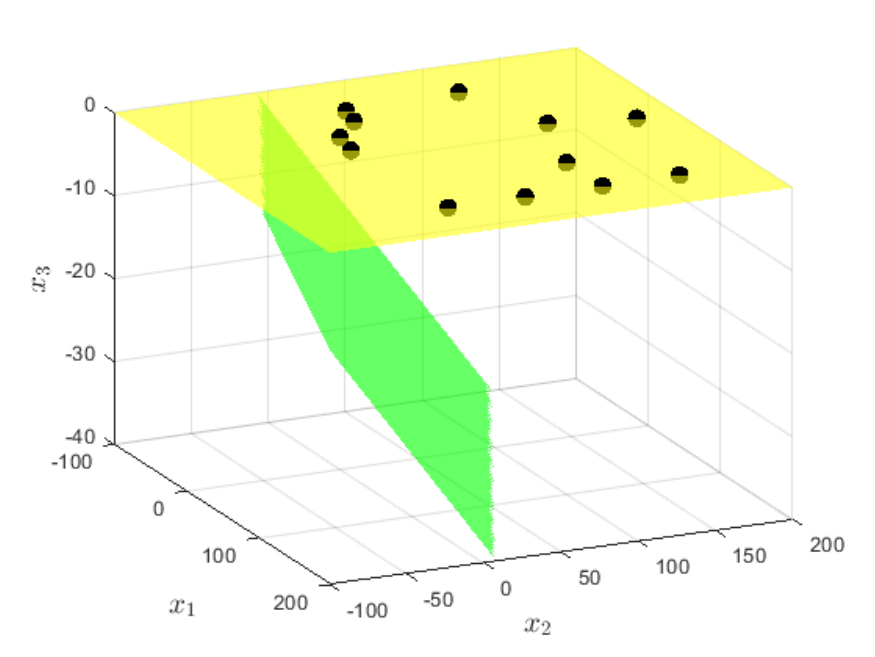}
			    \includegraphics[scale=.4]{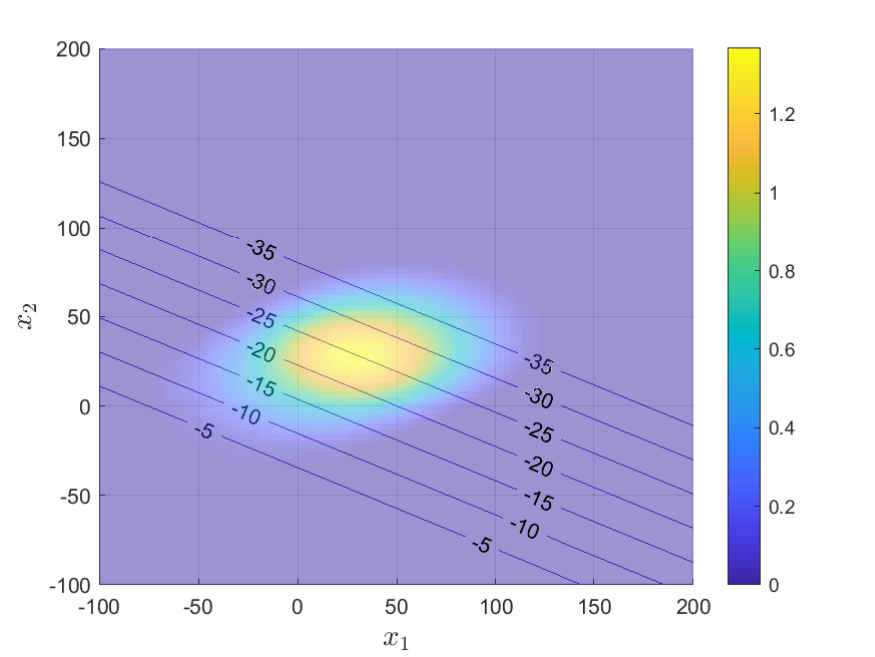}
    \caption{Left: the fault $\Gamma_{\tilde{\bm} } $
		in green and the surface measurement points $P_j$ 
		(sketched as black dots on the plane $x_3=0$ which is itself shown in yellow)
		in the case where $M=12$. Right: the slip $\cg$ on $\Gamma_{\tilde{\bm} }$
		viewed from above. $\cg$ 
		is taken to be in the direction of steepest descent, so only the magnitude is shown. Depth lines on $\Gamma_\bm$  are shown.}
    \label{forcing}
\end{figure}

\begin{figure}[htbp]
    \centering
      \includegraphics[scale=.4]{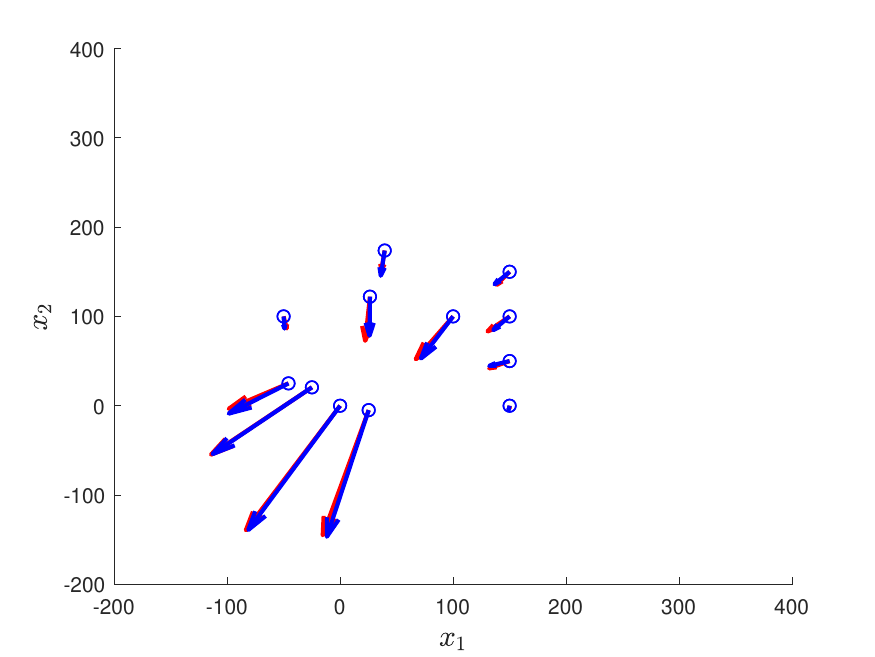}
			    \includegraphics[scale=.4]{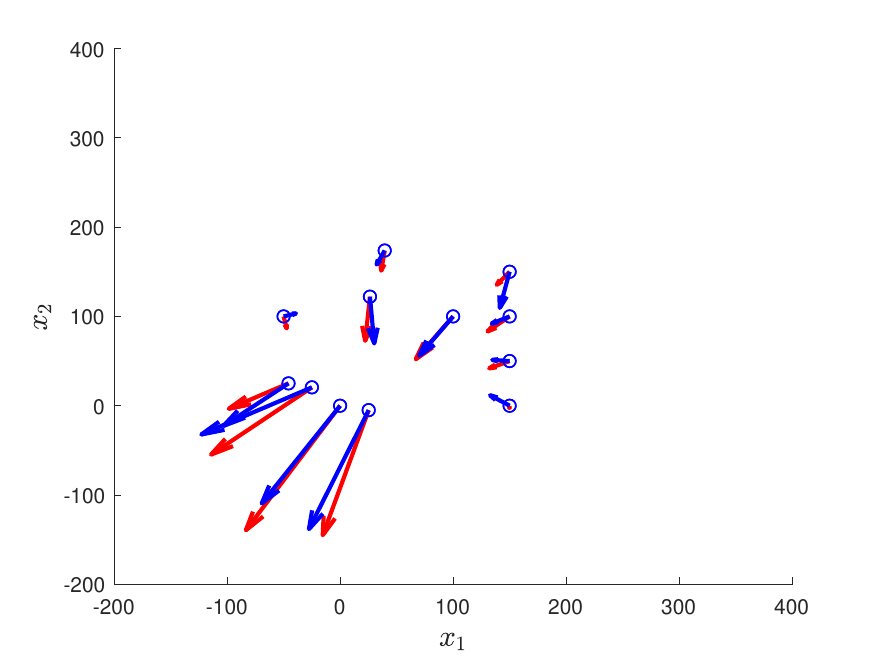}
      \includegraphics[scale=.4]{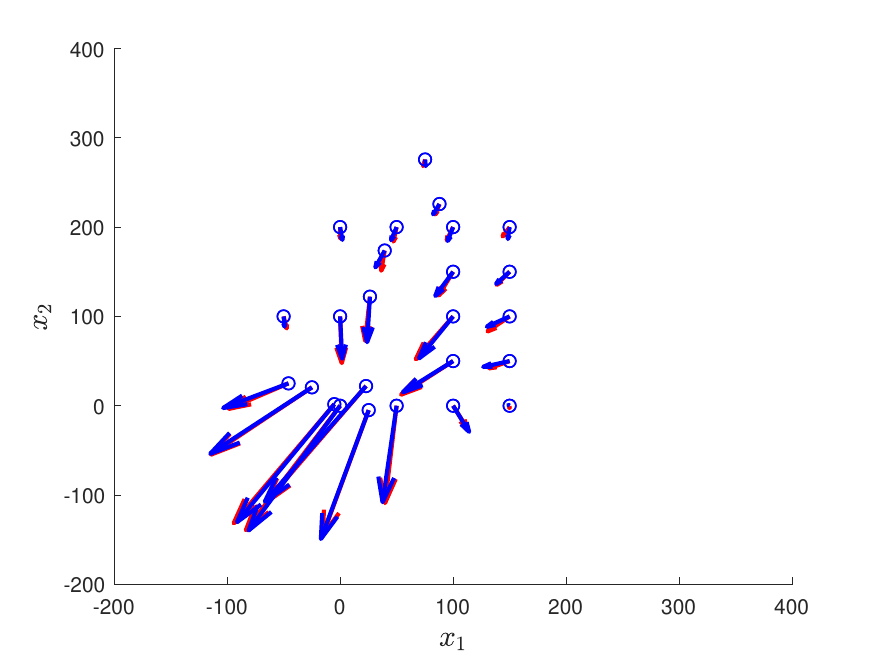}
      \includegraphics[scale=.4]{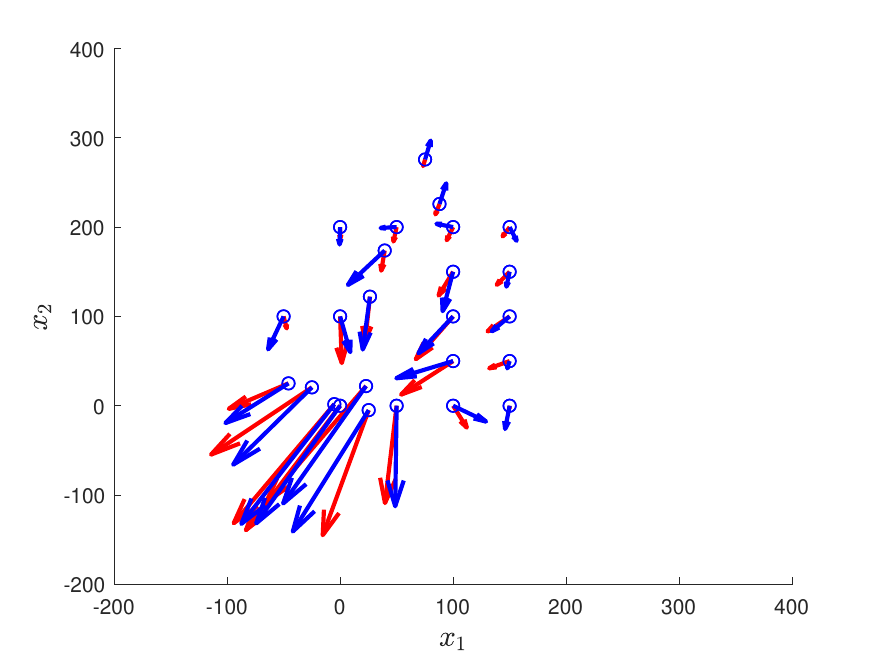}
			      \includegraphics[scale=.4]{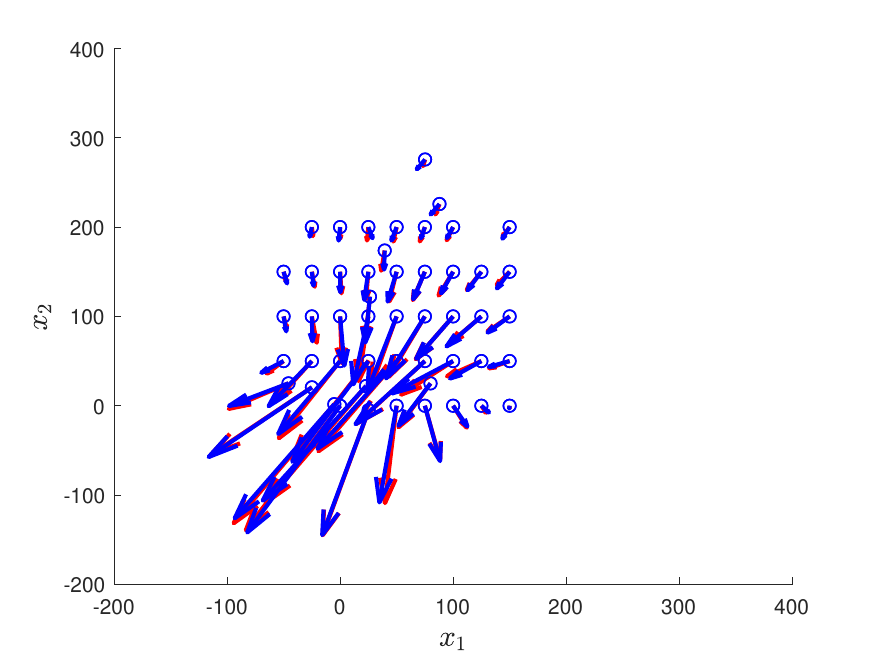}
	      \includegraphics[scale=.4]{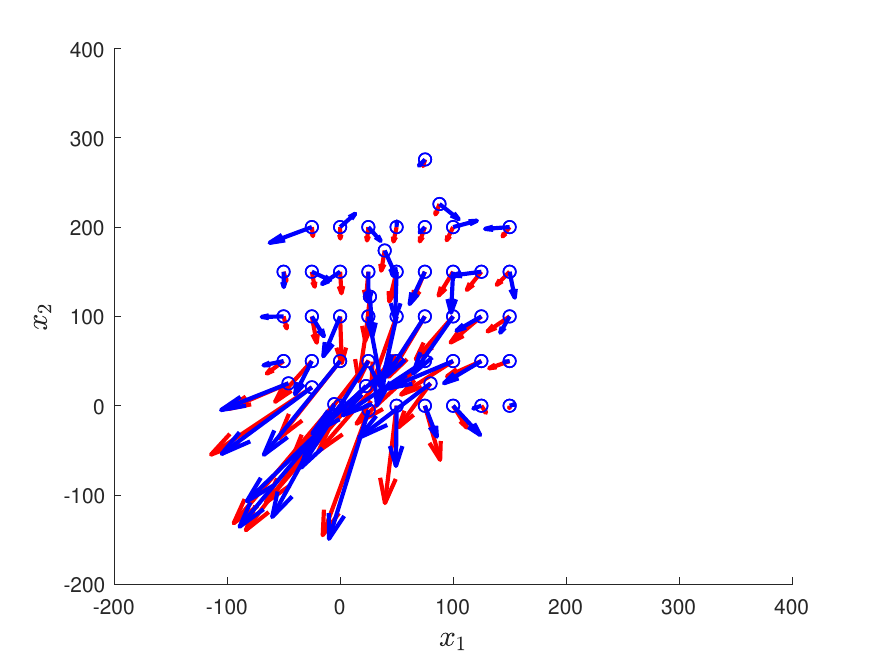}
    \caption{The data $\tilde{\bu}$ at the points $P_j$ (represented as circles) for the six inverse problem 
		solved numerically in this paper. Only the horizontal 
		displacements are sketched for clarity. Row 1 to 3: $M=12, 25, 50$. Column 1: low 
		$\sigma$, column 2: high $\sigma$.}
    \label{surface_disp}
\end{figure}

We consider data generated in a configuration
 closely related to  studies involving field data for a particular region and a specific seismic event 
\cite{volkov2019stochastic, volkov2017determining}.
To ensure that we perform a simulation
with realistic orders of magnitude,
the scaling is such that the unit for $\bx$ in $\RR^3$ is in kilometers, 
while $\tilde{\bu}$ and $\bg $ in  \eqref{noise model} 
are in meters, as in \cite{volkov2017determining}. 
For lower values of $M$ we will use a pattern of measurement points $P_j, j=1,.., M$
derived from the locations of in-situ measurement apparatus as 
set up by geophysicists \cite{volkov2017determining}. 
For generating forward data we pick the particular value
\bean
 \tilde{\bm} =
(   -0.12,  -0.26,  -14). \label{values}
\eean
We sketched $\Gamma_{\tilde{\bm} }$ in Figure \ref{forcing}, left,  where we also show the points 
$P_j$.    % in case $N=12$.
The slip field is assumed to be parallel to the steepest direction on $\Gamma_{\tilde{\bm} }$
and we assume that this is known when we solve the related inverse problem.
This direction of slip is characteristic of plate interface in subduction zones.
The slip field $\cg$ used for generating data as formulated by (\ref{noise model}-\ref{delta def})
%numerical simulations from this section 
is shown in
Figure \ref{forcing}, right. A very fine grid was used for computing the resulting 
displacement $\tilde{\bu}$.
The data for the inverse problem is the three dimensional displacements
 at the measurement points shown in 
Figure \ref{surface_disp}  to which we added
 Gaussian noise with covariance $\sigma^2 I$.
We present in this paper results for two values of $\sigma$, and for each value
of $\sigma$ three values of $(M, \dim H)$:
$(12, 27^2),  (25,37^2) ,  (50, 51^2) $.
In the lower  noise case, $\sigma$   was set to be equal to 5\%
of the maximum of the absolute values of the components of $\bu$
(in other words, 5\% of $\| \tilde{\bu} \|_{\infty}$).
For the particular realization used in solving the inverse problem, 
this led to a relative error in Euclidean norm of about 6\%.
In the higher noise case scenario $\sigma$ was set to be equal to 25\%
of the maximum of the absolute values of the components of $\bu$
(in other words, 25\% of $\| \tilde{\bu}\|_{\infty}$)
and  this led to 
 a relative error in Euclidean norm of about 30\%.
Both realizations are shown in Figure \ref{surface_disp}
(only the horizontal components are sketched for the sake of brevity).

\subsection{Numerical solution to the inverse problem} \label{Numerical solution}
\begin{figure}[htbp]
    \centering
      \includegraphics[scale=.4]{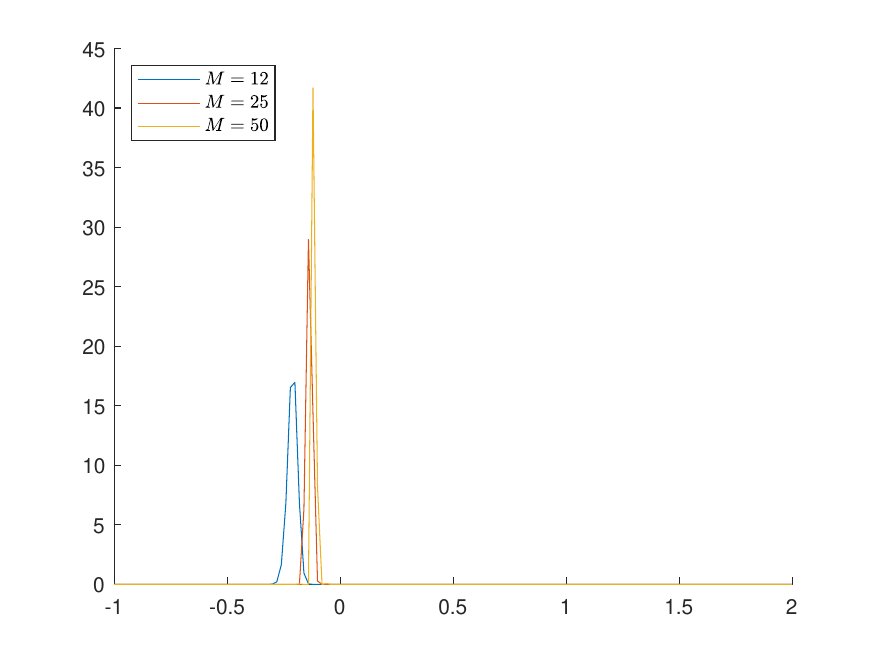}
			    \includegraphics[scale=.4]{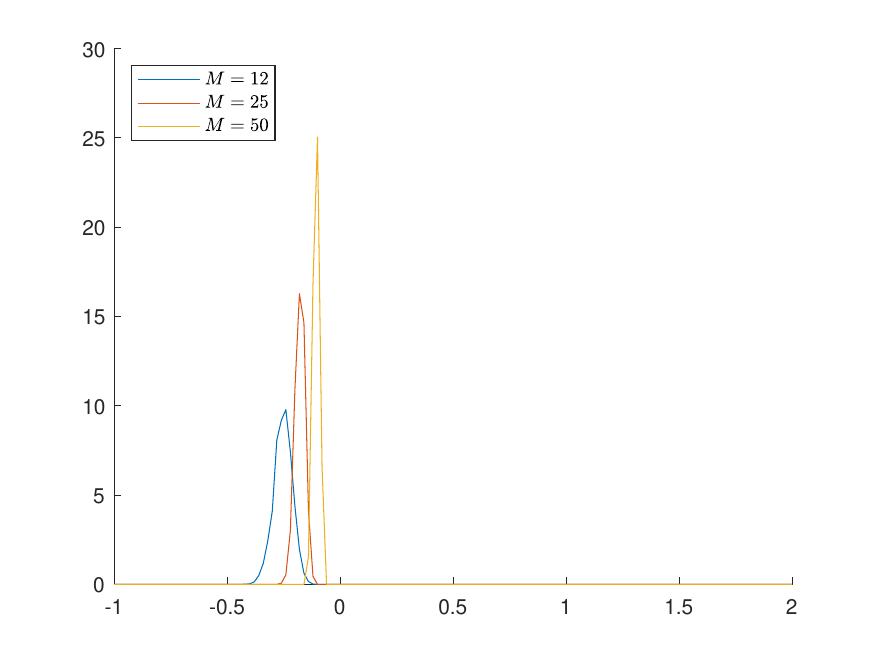}
      \includegraphics[scale=.4]{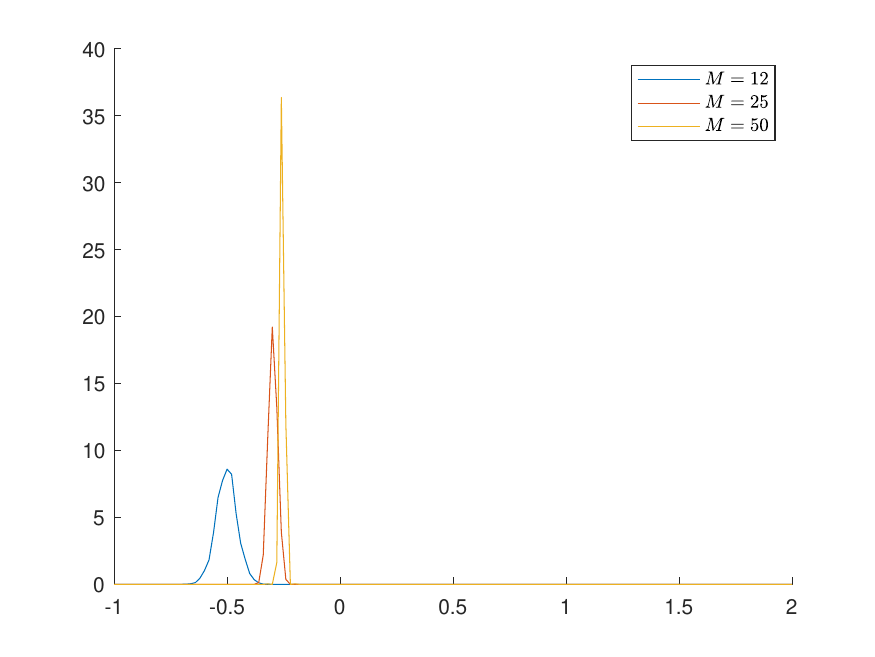}
      \includegraphics[scale=.4]{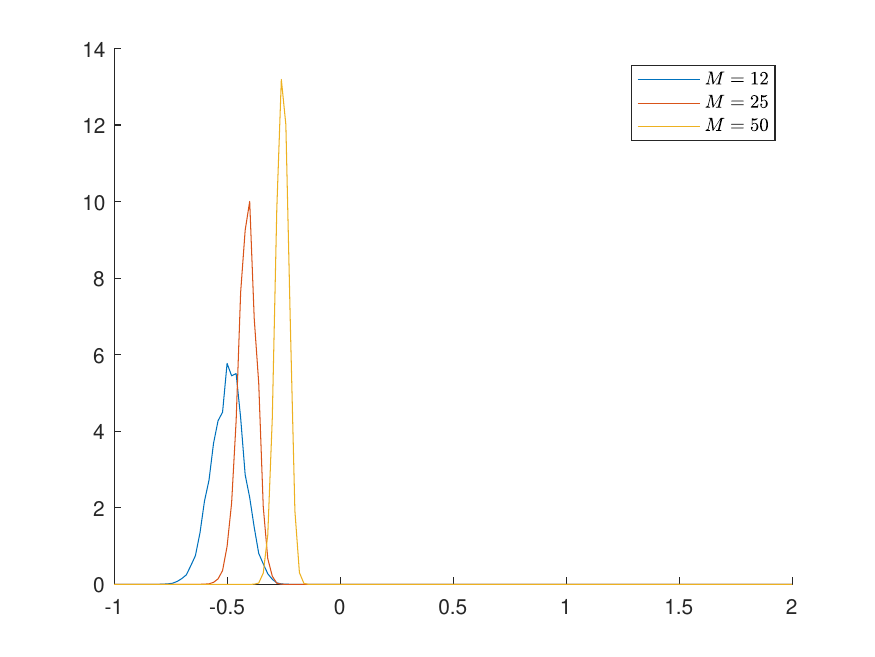}
			      \includegraphics[scale=.4]{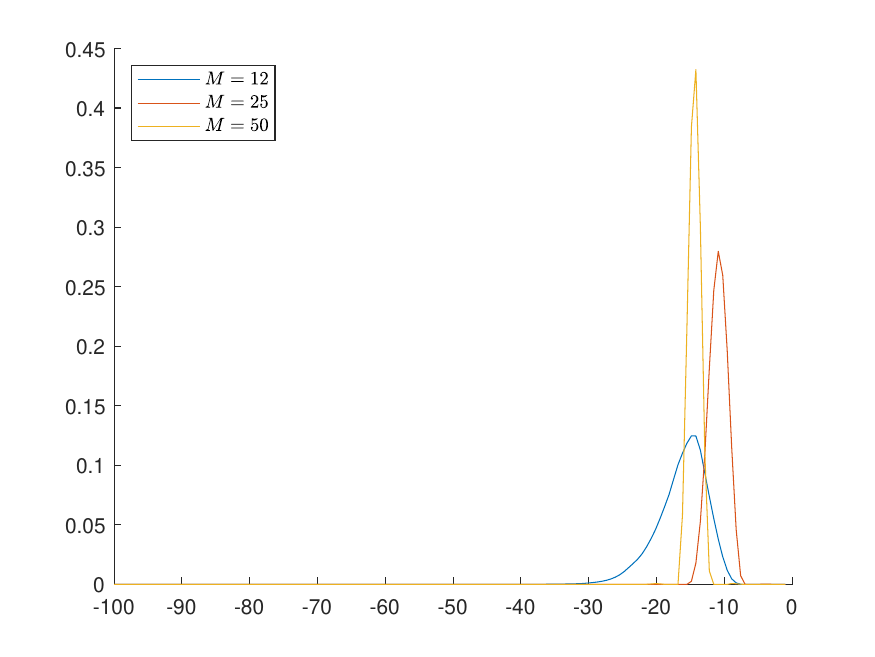}
	      \includegraphics[scale=.4]{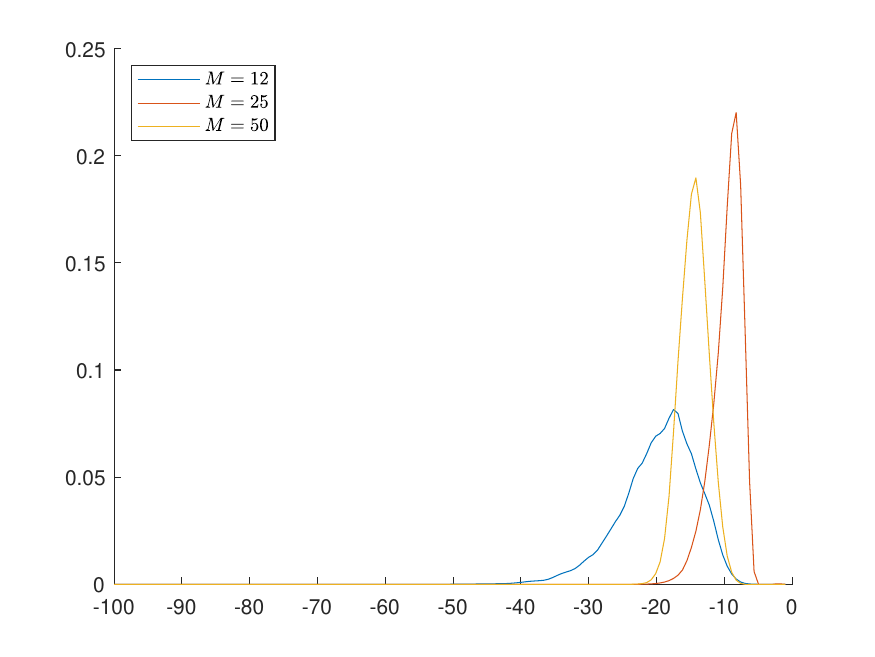}
    \caption{The  marginal posteriors of $a, b, d$ 
		computed from the data shown in Figure  \ref{surface_disp} in each of the six cases.
		Row 1, 2, 3: $a,b,d$.
		Left column: low noise case. Right column: high noise case.
In each figure, the computed posterior for $M=12, 25,$ and 50 are sketched.}
    \label{pdfs}
\end{figure}

\begin{figure}[htbp]
    \centering
      \includegraphics[scale=.4]{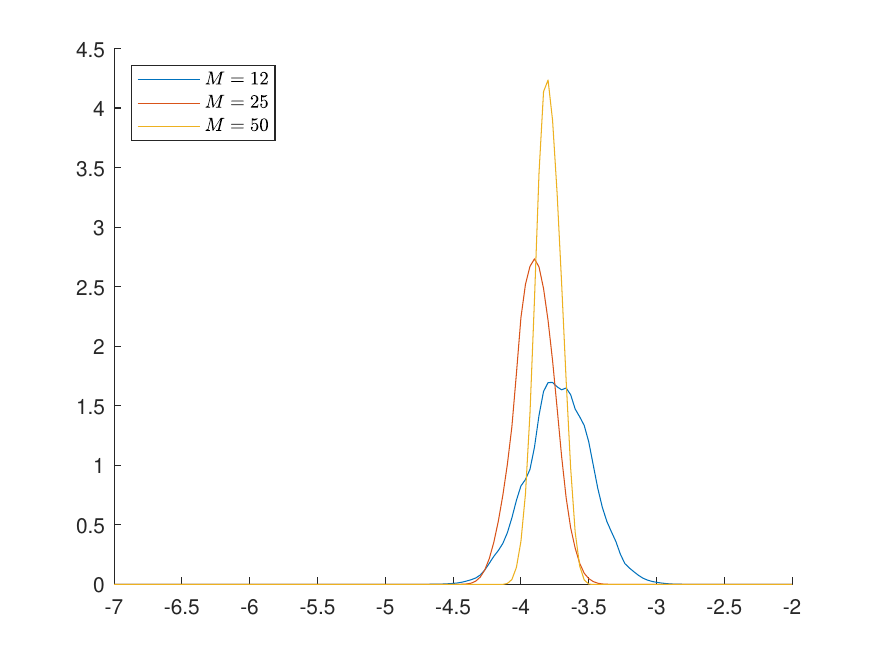}
			    \includegraphics[scale=.4]{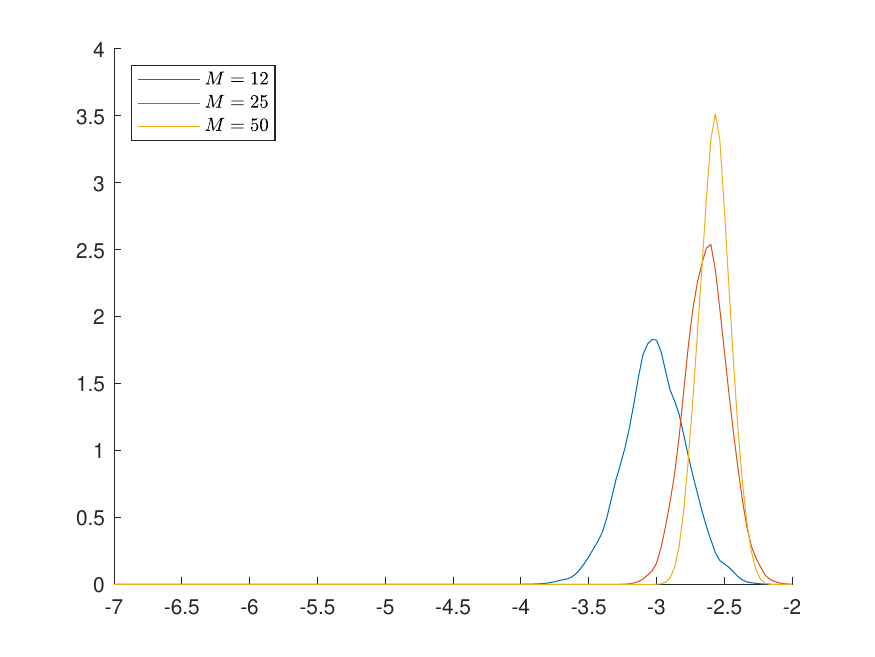}
        \caption{The computed marginal posteriors of the regularization parameter $C$. 
Left: low noise case. Right: high noise case.}
%In each figure, the computed posterior for $M_N=12, 25,$ and 50 are sketched. }
    \label{Cpdfs}
\end{figure}

We computed the posterior probability distribution function 
$\rho (\bm , C| \tilde{\bu}) $
given by (\ref{real1}-\ref{real2}) in each of the 6 cases introduced in the previous section.
The prior distribution of $\bm$ was taken to be uniform in the box
$[-1,2] \times [-1,2] \times [-100,-1]$, and the prior distribution of
$\log_{10} C$ was taken to be uniform in $[-7,-2]$. The computation of
 $\rho (\bm , C| \tilde{\bu}) $ was performed using the method of choice sampling,
more specifically, we used a modified version of the Metropolis algorithm
which is well suited to parallel computing \cite{calderhead2014general}.
In \cite{volkov2020stochastic}, section 4, we wrote explicitly a form of this algorithm for
computing  $\rho (\bm , C| \tilde{\bu}) $. Note that 
although the inverse problem considered 
in \cite{volkov2020stochastic} was similar to the one studied here, the nonlinear parameter to be reconstructed was different. \\
%In both studies, the novel idea was to use the posterior given by 
%$\rho (\bm , C| \tilde{\bu}) $: there is a great advantage to model $C$ as a random variable
%rather than to attempt to fix it using the maximum likelihood method, the discrepancy principle, or
%generalized cross validation.\\
The computed marginal posteriors of $\bm=(a,b,d)$ are graphed in Figure \ref{pdfs}.
This figure shows how these computed posteriors tighten around the  value
of $\tilde{\bm}$ as $M$
and $\dim H$ increase
as expected from Theorem \ref{main cv theorem}.
This tightening around the correct value occurs in both low and high $\sigma$
cases and is more narrow in the lower noise case.
Deciding that  the regularization parameter $C$  be a random variable
results in a more extensive random walk, however 
 it has the distinct advantage of 
sweeping through the entire range of values of $C$ 
in the support of the prior of $C$.
The computed marginal posteriors of $C$ are sketched in Figure
\ref{Cpdfs}. This figure illustrates how the algorithm automatically favors 
optimal values for $C$ depending on $M$, $\dim H$, and the noise level. % $\sigma$.

\subsection{Failure at fixed C} \label{Failure at fixed C}
\begin{figure}[htbp]
    \centering
      \includegraphics[scale=.4]{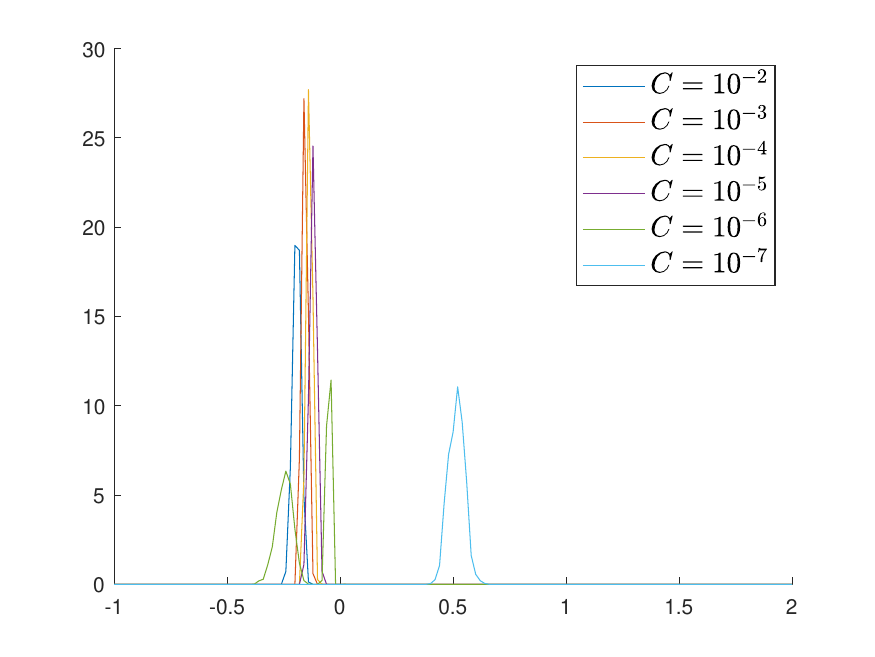}
			    \includegraphics[scale=.4]{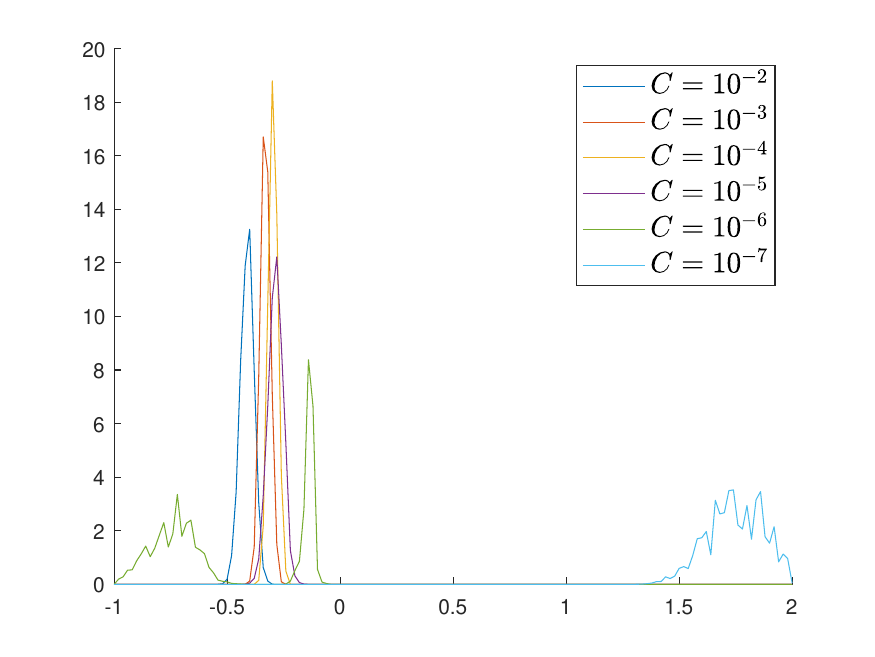}
    	    \includegraphics[scale=.4]{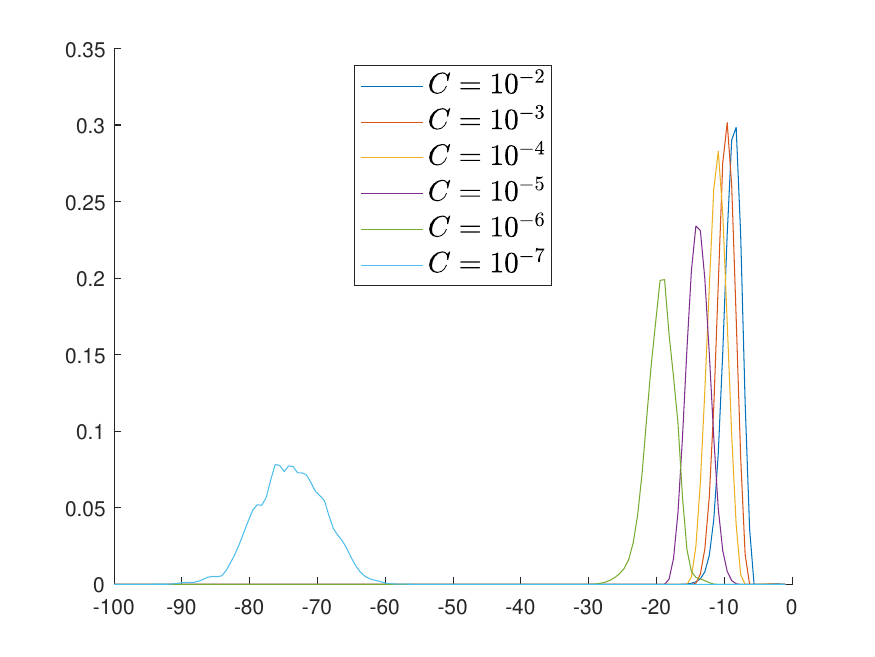}
        \caption{Computed posterior marginals assuming various fixed values of the regularization
				parameter $C$. First row: $a, b$. Second row: $d$. The computed marginals are only shown for the case $M=25$,  low $\sigma.$}
    \label{fixedC}
\end{figure}

% in C:\Users\darko\Desktop\RESEARCH\Projects for 2020\calc for third summer paper\25 points and low p
% run auto_fixed_C_august2020.m
% then
% post_process_for_fixedC.m
% and
% plot_results_fixedC.m
% a_figure_25_low
Fixing  a value for the regularization parameter $C$ is commonly done in linear inverse problems.
Often times, a value for $C$ is fixed in such a way that the solution displays satisfactory qualitative
features.   $\log_{10} C$ is varied linearly until such features appear.
Alternatively, one can use more objective criteria for selecting $C$ such as 
the maximum likelihood method, the discrepancy principle, or the 
generalized cross validation criterion.
However, the fault inverse problem is nonlinear in $\bm$. If one were to apply any of these 
methods to select a fixed $C$, the selection would depend on $\bm$ and as a result different 
candidates for $\bm$ would be unfairly compared, \cite{volkov2020stochastic}.
Better results are obtained if we fix the same value for $C$ for all $\bm$ in $B$, 
 \cite{volkov2019stochastic}.
Even then, determining  the optimal value for $C$ is not possible.
To illustrate this point, we plotted in Figure \ref{fixedC}
the computed posterior marginals assuming various fixed values of the regularization
				parameter $C$.
				Qualitatively, it appears that the values $10^{-7}$ and $10^{-6}$ have to be rejected,
				but it is unclear which of the remaining values is most suitable. 
				This example illustrates that modeling $C$ as a random variable and using
				a Bayesian approach built from the distribution function
			(\ref{real1}-\ref{real2})
				leads to far superior results.

%\newpage
%\section{Conclusion and perspectives for future work}
\vskip 10pt
\Large{\bf{Funding}} \\
\normalsize
This work was supported by
 Simons Foundation Collaboration Grant [351025].
%\vskip 10pt
\newpage
\bibliography{ref}{}
\bibliographystyle{abbrv}
%\bibliographystyle{apalike}
 
%\bibitem{S} E. P. Stephan, A boundary integral equation method for three-dimensional crack problems in elasticity, Mathematical Methods In The Applied Sciences, 1986 ,Volume: 8 Issue: 4.

\end{document}